# INVARIANT STATES AND RATES OF CONVERGENCE FOR A CRITICAL FLUID MODEL OF A PROCESSOR SHARING QUEUE

By Amber L. Puha[1] and Ruth J. Williams[2]

*California State University and University of California*

This paper contains an asymptotic analysis of a fluid model for a heavily loaded processor sharing queue. Specifically, we consider the behavior of solutions of critical fluid models as time approaches $\infty$. The main theorems of the paper provide sufficient conditions for a fluid model solution to converge to an invariant state and, under slightly more restrictive assumptions, provide a rate of convergence. These results are used in a related work by Gromoll for establishing a heavy traffic diffusion approximation for a processor sharing queue.

**1. Introduction.** This paper is a sequel to [10], which establishes a fluid (or functional law of large numbers) approximation for a heavily loaded processor sharing queue. In [10], a stochastic process $\mu(\cdot)$ taking values in $\mathcal{M}_F$, the space of finite, nonnegative Borel measures on $\mathbb{R}_+ = [0, \infty)$ endowed with the topology of weak convergence, is used to track the evolution in time of the state of a processor sharing queue. At time $t$, $\mu(t)$ is the measure that has one unit of mass at the residual service time of each job present in the system at time $t$. From the *measure-valued state descriptor* $\mu(\cdot)$, one can recover the traditional performance processes, such as the queue length and workload processes (cf. [10], Section 2.3). Under mild conditions, it is proved in [10] that the fluid scaled state descriptors for a sequence of heavily loaded processor sharing queues converge in distribution to a measure-valued stochastic process, which we refer to as a *fluid limit* (cf. [10], Theorem 3.2). Almost every sample path of this fluid limit is a solution of a certain (deterministic)

Received September 2002; revised April 2003.
[1]Supported in part by a University of California Office of the President Postdoctoral Fellowship and an NSF Mathematical Sciences Postdoctoral Fellowship.
[2]Supported in part by NSF Grant DMS-00-71408, a John Simon Guggenheim Fellowship and a gift from the David and Holly Mendel fund.

*AMS 2000 subject classifications.* Primary 60K25; secondary 68M20, 90B22.

*Key words and phrases.* Processor sharing queue, critical fluid model, measure-valued solution, invariant states, coupling renewal processes, renewal functions, renewal measures.

critical fluid model. In this paper, we study the asymptotic behavior as time tends to $\infty$ of the solutions of this critical fluid model.

In [1] and [2], Bramson studied the asymptotic behavior of solutions of critical fluid models associated with open multiclass queueing networks operating under two HL (head-of-the-line) service disciplines. Then, in [3], Bramson showed that, if the critical fluid model associated with an open multiclass HL queueing network has a certain asymptotic property, then a condition known as state space collapse holds. In a companion work to [3], Williams [14] showed that state space collapse, plus an algebraic condition on the first-order queueing model data, is sufficient to imply a heavy traffic diffusion approximation for an open multiclass queueing network operating under an HL service discipline. To illustrate this modular approach, Bramson [3] and Williams [14] applied their results, together with the results of [1] and [2], to obtain new heavy traffic diffusion limit theorems for FIFO networks of Kelly type and for networks with an HLPPS (head-of-the-line proportional processor sharing) service discipline. Processor sharing, as considered in this paper, is not an HL service discipline. However, an analogue of the modular approach of [3] and [14] is developed for a processor sharing queue in [9]. The results proved here are used in [9] to prove a state space collapse result, which in turn is used in [9] to establish a heavy traffic diffusion approximation for a processor sharing queue.

To state our results, we need to recall the description of the *critical fluid model* from [10]. The model has two parameters, $\alpha \in (0, \infty)$ and a Borel probability measure $\nu$ on $\mathbb{R}_+$ that does not charge the origin [$\nu(\{0\}) = 0$] and has a finite first moment [$\int_{\mathbb{R}_+} x\nu(dx) < \infty$]. These parameters correspond to parameters in the queueing system. Specifically, $\alpha$ corresponds to the long-run average rate at which jobs arrive to the system, and the probability measure $\nu$ corresponds to the distribution of the i.i.d. service times for those jobs. The qualifier *critical* refers to the fact that we are interested in the critically loaded regime where the service and arrival rates are equal. Thus, it is assumed throughout that

$$(1.1) \qquad \alpha = \left(\int_{\mathbb{R}_+} x\nu(dx)\right)^{-1}.$$

The pair $(\alpha, \nu)$ is referred to as the data for the critical fluid model, or simply the *critical data*. Here we only consider solutions of a critical fluid model, and we simply refer to these as fluid model solutions. In particular, the assumption of critical data is implicit.

A fluid model solution $\bar{\mu}(\cdot)$ is a deterministic function of time, taking values in $\mathcal{M}_F$, that satisfies conditions (C1)–(C4). To state these conditions, we need to introduce some notation. For a Borel set $A \subset \mathbb{R}_+$, let $\mathbb{1}_A$ denote the indicator function of the set $A$. To simplify the notation, we use the shorthand notation $\mathbb{1}$ in place of $\mathbb{1}_{\mathbb{R}_+}$. For $\zeta \in \mathcal{M}_F$, the real-valued projection of $\zeta$ associated with a bounded, real-valued, Borel measurable function $g$ defined on $\mathbb{R}_+$ is denoted by $\langle g, \zeta \rangle = \int_{\mathbb{R}_+} g(x)\zeta(dx)$. The dynamic conditions [see (C3)] that an $\mathcal{M}_F$-valued function $\bar{\mu}(\cdot)$ must satisfy in order to be a fluid model solution involve the real-valued projections of $\bar{\mu}(\cdot)$ over the class of functions

$$\mathcal{C} = \{g \in \mathbf{C}_b^1(\mathbb{R}_+) : g(0) = 0, g'(0) = 0\}.$$

Here $\mathbf{C}_b^1(\mathbb{R}_+)$ denotes the space of once continuously differentiable real-valued functions defined on $\mathbb{R}_+$ that, together with their first derivatives, are bounded on $\mathbb{R}_+$. The requirement that $g$ and $g'$ vanish at the origin is imposed to avoid possible singular behavior of $\langle g, \bar{\mu}(\cdot) \rangle$ and $\langle g', \bar{\mu}(\cdot) \rangle$, associated with mass in the fluid model abruptly disappearing as it reaches the origin. Such behavior corresponds to jobs in the queueing system abruptly departing when their residual service times reach 0.

A *fluid model solution* is a function $\bar{\mu} : [0, \infty) \longrightarrow \mathcal{M}_F$ that satisfies the following four conditions.

(C1) The function $\bar{\mu}(\cdot)$ is continuous.
(C2) For each $t \geq 0$, $\langle \mathbb{1}_{\{0\}}, \bar{\mu}(t) \rangle = 0$.
(C3) For each $g \in \mathcal{C}$, $\bar{\mu}(\cdot)$ satisfies

$$\langle g, \bar{\mu}(t) \rangle = \langle g, \bar{\mu}(0) \rangle - \int_0^t \frac{\langle g', \bar{\mu}(s) \rangle}{\langle \mathbb{1}, \bar{\mu}(s) \rangle} \, ds + \alpha t \langle g, \nu \rangle \tag{1.2}$$

for all $0 \leq t < t^* = \inf\{s \geq 0 : \langle \mathbb{1}, \bar{\mu}(s) \rangle = 0\}$.
(C4) For all $t \geq t^*$, $\langle \mathbb{1}, \bar{\mu}(t) \rangle = 0$.

See [10], Section 3.1, for an interpretation of (C1)–(C4) in terms of the dynamics of a processor sharing queue. In fact, using dominated convergence and (C2), it is straightforward to see that $\bar{\mu} : [0, \infty) \longrightarrow \mathcal{M}_F$ satisfies (C1)–(C4) if and only if it satisfies these conditions with $\mathcal{C}$ replaced by $\tilde{\mathcal{C}} = \{g \in \mathbf{C}_b^1(\mathbb{R}_+) : g(0) = 0\}$. The more restictive class was used in [10] as it simplified the proof of the existence of solutions. In addition, as is proved in [10] and explained below, for the nontrivial fluid model solutions considered here, $t^* = \infty$.

To facilitate the present discussion, we review some results from [10] concerning fluid model solutions. Let

$$\mathcal{M}_F^c = \{\xi \in \mathcal{M}_F : \langle \mathbb{1}_{\{x\}}, \xi \rangle = 0 \text{ for all } x \in \mathbb{R}_+\},$$

where $c$ stands for continuous. Theorem 3.1 in [10] states that, for each measure $\xi \in \mathcal{M}_F^c$, there exists a unique fluid model solution $\bar{\mu}_\xi(\cdot)$ such that $\bar{\mu}_\xi(0) = \xi$. If $\xi = \mathbf{0}$, where $\mathbf{0}$ denotes the zero measure, then, by (C4), $\bar{\mu}_\xi(\cdot) \equiv \mathbf{0}$. Let

$$\mathcal{M}_F^{c,p} = \{\xi \in \mathcal{M}_F^c : \xi \neq \mathbf{0}\},$$

where $p$ stands for positive. In [10], it was also shown that, for $\xi \in \mathcal{M}_{\mathrm{F}}^{c,p}$, $\bar{\mu}_\xi(t) \in \mathcal{M}_{\mathrm{F}}^{c,p}$ for all $t \geq 0$ (cf. [10], Theorem 3.1 and Proposition 4.6). In particular, if $\xi \in \mathcal{M}_{\mathrm{F}}^{c,p}$, then $t^* = \infty$.

Given $\xi \in \mathcal{M}_{\mathrm{F}}^c$, it is natural to ask about the asymptotic behavior of $\bar{\mu}_\xi(t)$ as $t$ tends to $\infty$. Specifically, as $t$ tends to $\infty$, does $\bar{\mu}_\xi(t)$ converge in some sense? If so, what is the limit and how fast is the convergence? To answer these questions, we begin by identifying the possible limiting measures. Extending the terminology in [3] to the present setting, a measure $\xi \in \mathcal{M}_{\mathrm{F}}^c$ is said to be an *invariant state* if

$$\bar{\mu}_\xi(t) = \xi \qquad \text{for all } t \geq 0.$$

Similarly, the collection of invariant states $\mathfrak{I}$, which is given by

$$\mathfrak{I} = \{\xi \in \mathcal{M}_{\mathrm{F}}^c : \bar{\mu}_\xi(t) = \xi \text{ for all } t \geq 0\},$$

is called the *invariant manifold*. Here it turns out that the invariant manifold $\mathfrak{I}$ is a one-parameter family of measures that is determined by the probability measure $\nu$. To describe $\mathfrak{I}$, we need to introduce some notation. Let $F$ denote the cumulative distribution function associated with the probability measure $\nu$. The distribution function $F$ has associated with it an excess lifetime cumulative distribution function $F_{\mathrm{e}}$, which is given by

$$F_{\mathrm{e}}(x) = \alpha \int_0^x (1 - F(y)) \, dy \qquad \text{for all } x \in \mathbb{R}_+.$$

In particular, $F_{\mathrm{e}}$ has probability density function

$$f_{\mathrm{e}}(x) = \alpha(1 - F(x)) \qquad \text{for all } x \in \mathbb{R}_+.$$

Note that (1.1) was used to simplify the form of the normalizing constant here. Let $\nu_{\mathrm{e}}$ denote the Borel probability measure on $\mathbb{R}_+$ that has density function $f_{\mathrm{e}}$, that is, $\langle \mathbb{1}_{[0,x]}, \nu_{\mathrm{e}} \rangle = F_{\mathrm{e}}(x) = \int_0^x f_{\mathrm{e}}(y) \, dy$ for all $x \in \mathbb{R}_+$. We call $\nu_{\mathrm{e}}$ the *excess lifetime probability measure*. Define

$$\beta_{\mathrm{e}} = \frac{1}{\langle \chi, \nu_{\mathrm{e}} \rangle},$$

where $\chi(x) = x$ for $x \in \mathbb{R}_+$. The right member above is interpreted as 0 if the first moment of $\nu_{\mathrm{e}}$ is infinite.

THEOREM 1.1. *A measure $\xi \in \mathcal{M}_{\mathrm{F}}^c$ is an invariant state if and only if $\xi = c\nu_{\mathrm{e}}$ for some $c \in [0, \infty)$. Equivalently, the invariant manifold $\mathfrak{I}$ is given by*

$$\mathfrak{I} = \{c\nu_{\mathrm{e}} : c \in [0, \infty)\}.$$

Theorem 1.1 is proved in Section 3.

Let $\xi \in \mathcal{M}_\text{F}^c$. We wish to identify conditions under which $\bar{\mu}_\xi(t)$ converges to a point on the invariant manifold as $t$ tends to $\infty$ and to determine the limiting state. For this, we define the fluid analogue of the workload at time $t \in [0, \infty)$ to be given by $\langle \chi, \bar{\mu}_\xi(t) \rangle$. By Theorem 3.1 in [10], $\langle \chi, \bar{\mu}_\xi(t) \rangle = \langle \chi, \xi \rangle$ for all $t \geq 0$. [This holds even if $\langle \chi, \xi \rangle = \infty$, in which case $\langle \chi, \bar{\mu}_\xi(t) \rangle = \infty$ for all $t \geq 0$.] Thus, when $\bar{\mu}_\xi(t)$ converges to an element $c\nu_\text{e}$ in $\mathcal{I}$ as $t$ tends to $\infty$ and both $\langle \chi, \xi \rangle$ and $\langle \chi, \nu_\text{e} \rangle$ are finite, one might expect the first moment, $c\langle \chi, \nu_\text{e} \rangle$, of the limit to be given by $\langle \chi, \xi \rangle$, or, equivalently, that $c = \beta_\text{e} \langle \chi, \xi \rangle$. Indeed, we have the following result.

THEOREM 1.2. *Let $\xi \in \mathcal{M}_\text{F}^c$. If $\langle \chi, \xi \rangle < \infty$, then $\bar{\mu}_\xi(t)$ converges weakly to $\beta_\text{e} \langle \chi, \xi \rangle \nu_\text{e}$ as $t \to \infty$.*

Notice that $\beta_\text{e} > 0$ if and only if $\langle \chi^2, \nu \rangle < \infty$, since, for $\gamma \in \mathbb{R}_+$, $\langle \chi^\gamma, \nu_\text{e} \rangle < \infty$ if and only if $\langle \chi^{\gamma+1}, \nu \rangle < \infty$. Therefore, the case in which $\xi \neq \mathbf{0}$ and $\langle \chi^2, \nu \rangle = \infty$ is degenerate in the sense that $\bar{\mu}_\xi(t)$ converges to the zero measure as $t$ tends to $\infty$, but $\langle \chi, \bar{\mu}_\xi(t) \rangle$ does not converge to 0.

The result in Theorem 1.2 is more general than, but consistent with, Proposition 5 of [4], which concerns the asymptotic behavior of a fluid approximation for the queue length of a heavily loaded processor sharing queue. Theorem 1.2 is proved in Section 4, using proof techniques similar to those employed in [4].

Finally, we wish to give a rate at which $\bar{\mu}_\xi(t)$ converges as $t$ tends to $\infty$. In fact, we will prove two rate of convergence results. The first gives a rate of convergence in terms of a metric on $\mathcal{M}_\text{F}$ that induces the weak topology. For this, let $\rho$ denote the extension of the Prohorov metric to $\mathcal{M}_\text{F}$. Specifically, for $\zeta_1, \zeta_2 \in \mathcal{M}_\text{F}$, $\rho(\zeta_1, \zeta_2)$ is given by

$$\rho(\zeta_1, \zeta_2) = \inf\{\delta > 0 : \langle \mathbb{1}_B, \zeta_1 \rangle \leq \langle \mathbb{1}_{B^\delta}, \zeta_2 \rangle + \delta$$
(1.3)
$$\text{and } \langle \mathbb{1}_B, \zeta_2 \rangle \leq \langle \mathbb{1}_{B^\delta}, \zeta_1 \rangle + \delta,$$
$$\text{for all nonempty, closed sets } B \subset \mathbb{R}_+\},$$

where, for each nonempty, closed set $B \subset \mathbb{R}_+$,

$$B^\delta = \left\{ x \in \mathbb{R}_+ : \inf_{y \in B} |x - y| < \delta \right\}.$$

Note that, under $\rho$, $\mathcal{M}_\text{F}$ is a Polish space. Moreover, if $\{\zeta_n, n = 1, 2, \ldots\} \subset \mathcal{M}_\text{F}$ and $\zeta \in \mathcal{M}_\text{F}$, then $\zeta_n$ converges weakly to $\zeta$ as $n$ tends to $\infty$ if and only if $\lim_{n \to \infty} \rho(\zeta_n, \zeta) = 0$ (cf. [6], Chapter 3, Theorems 1.7 and 3.1, which readily generalize from the set of Borel probability measures to $\mathcal{M}_\text{F}$). Our second rate result gives a rate of convergence in terms of the total variation distance. For a signed, Borel measure $\zeta$ on $\mathbb{R}_+$,

$$\begin{aligned}(1.4)\quad \|\zeta\|_{\text{TV}} = \sup\{|\langle g,\zeta\rangle| \text{ such that } &g:\mathbb{R}_+ \longrightarrow \mathbb{R} \text{ is Borel measurable}\\ &\text{and } |g(x)| \leq 1 \text{ for all } x \in \mathbb{R}_+\}.\end{aligned}$$

Note that, if $\{\zeta_n, n=1,2,\dots\} \subset \mathcal{M}_{\text{F}}$, $\zeta \in \mathcal{M}_{\text{F}}$ and $\lim_{n\to\infty} \|\zeta_n - \zeta\|_{\text{TV}} = 0$, then $\zeta_n$ converges weakly to $\zeta$ as $n \to \infty$. However, the converse is not true in general (cf. [5], page 69). For each of our rate of convergence results, the convergence is uniform over sets of initial conditions satisfying certain moment constraints. These sets take the following form. For any finite, positive constants $\varepsilon$ and $M$, let

$$(1.5)\quad \mathcal{B}_\rho^{M,\varepsilon} = \{\xi \in \mathcal{M}_{\text{F}}^c : \langle \mathbb{1},\xi\rangle \vee \langle \chi,\xi\rangle \vee \langle \chi^{1+\varepsilon},\xi\rangle \leq M\},$$

$$(1.6)\quad \mathcal{B}_{\text{TV}}^{M,\varepsilon} = \{\xi \in \mathcal{M}_{\text{F}}^c : \langle \mathbb{1},\xi\rangle \vee \langle \chi,\xi\rangle \vee \langle \chi^2,\xi\rangle \vee \langle \chi^{2+\varepsilon},\xi\rangle \leq M\}.$$

Of course, for $\xi \in \mathcal{M}_F^c$, if $\langle \mathbb{1},\xi\rangle \vee \langle \chi^{1+\varepsilon},\xi\rangle \leq M$, then $\xi \in \mathcal{B}_\rho^{2M,\varepsilon}$. Similarly, if $\langle \mathbb{1},\xi\rangle \vee \langle \chi^{2+\varepsilon},\xi\rangle \leq M$, then $\xi \in \mathcal{B}_{\text{TV}}^{2M,\varepsilon}$. Definitions (1.5) and (1.6) are used to simplify the tracking of constants in our proofs. Our rate of convergence results are summarized by the following.

THEOREM 1.3. *Let $M$ be a fixed, finite, positive constant.*

(i) *If, for some $\varepsilon > 0$, $\langle \chi^{2+\varepsilon},\nu\rangle < \infty$, then there exist a finite, positive constant $C_\rho$ and a finite, positive time $T_\rho$ such that*

$$(1.7)\quad \sup_{\xi \in \mathcal{B}_\rho^{M,\varepsilon}} (\bar{\mu}_\xi(t), \beta_{\text{e}}\langle \chi,\xi\rangle \nu_{\text{e}}) \leq C_\rho t^{-\varepsilon/4} \quad \text{for all } t \geq T_\rho.$$

(ii) *If, for some $\varepsilon > 0$, $\langle \chi^{3+\varepsilon},\nu\rangle < \infty$, then there exist a finite, positive constant $C_{\text{TV}}$ and a finite, positive time $T_{\text{TV}}$ such that*

$$(1.8)\quad \sup_{\xi \in \mathcal{B}_{\text{TV}}^{M,\varepsilon}} \|\bar{\mu}_\xi(t) - \beta_{\text{e}}\langle \chi,\xi\rangle \nu_{\text{e}}\|_{\text{TV}} \leq C_{\text{TV}} t^{-\varepsilon} \quad \text{for all } t \geq T_{\text{TV}}.$$

In Theorem 1.3, the times $T_\rho$ and $T_{\text{TV}}$ and the constants $C_\rho$ and $C_{\text{TV}}$ depend on the values of the constants $M$ and $\varepsilon$ and on the critical data $(\alpha,\nu)$. In the proofs, we have not tried to obtain the best possible estimates for these constants.

In [9], Theorem 1.3(i) is used to prove a state space collapse result. In that application, it is the uniform convergence over sets of the form $\mathcal{B}_\rho^{M,\varepsilon}$ for $M, \varepsilon \in (0,\infty)$ that is critical. In fact, the specific rate and value of the constants are not important for the argument. Although Theorem 1.3(ii) is not needed for [9], we have included it here for its intrinsic interest and potential use in other applications.

The proof of part (i) of Theorem 1.3 exploits the asymptotic behavior of the renewal function for a zero-delayed renewal process with interarrival distribution determined by the probability measure $\nu_e$. The condition $\langle \chi^{2+\varepsilon}, \nu \rangle < \infty$ is slightly stronger than requiring that this interarrival distribution have a finite mean, that is, that $\beta_e > 0$. This condition is used in the proof to obtain a rate of convergence for Blackwell's renewal theorem. Similarly, the proof of part (ii) of Theorem 1.3 exploits the asymptotic behavior of the renewal measures for certain delayed renewal processes with interarrival distribution determined by the probability measure $\nu_e$. The condition $\langle \chi^{3+\varepsilon}, \nu \rangle < \infty$ is slightly stronger than requiring that the interarrival distribution have a finite second moment. This condition is used in the proof to obtain a rate at which the renewal measures converge in the total variation distance to the stationary renewal measure. Both the rate of convergence for Blackwell's renewal theorem and the rates of convergence for renewal measures rely on the coupling results developed in [12]. In using those results, we pay careful attention to the dependence of the various constants on the initial measure $\xi$ and the interarrival distribution $\nu_e$.

The remainder of the paper is devoted to the proofs of Theorems 1.1–1.3. Section 2 contains some background and two preparatory lemmas. Then Theorems 1.1, 1.2, 1.3(i) and 1.3(ii) are proved in Sections 3, 4, 5 and 6, respectively. In the Appendix, coupling results from [12] are applied to verify some of the estimates used in Section 6 to prove Theorem 1.3(ii).

**2. Background.** Recall that, for $\xi \in \mathcal{M}_F^c$, $\bar{\mu}_\xi(\cdot)$ denotes the unique fluid model solution such that $\bar{\mu}_\xi(0) = \xi$. Given $\xi \in \mathcal{M}_F^c$, the fluid analogue of the queue length $\bar{Z}(\cdot)$ is defined by

$$\bar{Z}(t) = \langle \mathbb{1}, \bar{\mu}_\xi(t) \rangle \qquad \text{for all } t \geq 0. \tag{2.1}$$

For obvious reasons, $\bar{Z}(t)$ is referred to as the *total mass* at time $t$. Due to (C1), $\bar{Z}(\cdot)$ is continuous. As previously noted, if $\xi \neq \mathbf{0}$, then $\bar{\mu}_\xi(t) \neq \mathbf{0}$ for all $t \geq 0$ (cf. [10], Theorem 3.1), and so $\bar{Z}(t)$ is strictly positive for all $t \geq 0$. Conversely, if $\xi = \mathbf{0}$, then $\bar{\mu}_\xi(\cdot) \equiv \mathbf{0}$ and $\bar{Z}(\cdot) \equiv 0$. Given this, for each $t \geq 0$, the fluid analogue of the cumulative service per job $\bar{S}(t)$ is defined by

$$\bar{S}(t) = \begin{cases} 0, & \text{if } \xi = \mathbf{0}, \\ \int_0^t (\bar{Z}(s))^{-1} ds, & \text{otherwise.} \end{cases} \tag{2.2}$$

Thus, at time $t \geq 0$, $\bar{S}(t)$ denotes the *cumulative service per unit of mass* in the system up to time $t$. Since $\bar{Z}(\cdot)$ is continuous and $\bar{Z}(t) > 0$ for all $t \geq 0$ when $\xi \neq \mathbf{0}$, it follows that $\bar{S}(\cdot)$ is continuously differentiable. The reader will note that, in order to avoid cluttering the notation, we choose not to append a subscript $\xi$ to quantities defined by (2.1) and (2.2), since it is typically clear from the context which fluid model solution is under consideration. In [10], it was shown that, if $\xi \in \mathcal{M}_{\mathrm{F}}^{c,p}$, then, for each $t \geq 0$ and $x \in \mathbb{R}_+$,

$$(2.3) \qquad \langle \mathbb{1}_{[0,x]}, \bar{\mu}_\xi(t) \rangle = \langle \mathbb{1}_{(\bar{S}(t), \bar{S}(t)+x]}, \xi \rangle + \int_0^t G^x(\bar{S}(t) - \bar{S}(s))\, ds,$$

where, for each $x \in \mathbb{R}_+$,

$$G^x(y) = f_{\mathrm{e}}(y) - f_{\mathrm{e}}(x+y) \qquad \text{for all } y \in \mathbb{R}_+$$

(cf. [10], Lemma 4.3 and (4.33)). Here we have used (C2). For each $t \geq 0$, this gives an explicit description of the measure $\bar{\mu}_\xi(t)$ in terms of the nonzero initial measure $\xi$ and the cumulative service per unit of mass function $\bar{S}(\cdot)$.

To state what is known about $\bar{S}(\cdot)$ for a given $\xi \in \mathcal{M}_{\mathrm{F}}^{c,p}$, we need to introduce the *renewal function* $U_{\mathrm{e}}(\cdot)$ associated with the critical data $(\alpha, \nu)$ and the *truncated initial workload function* $H_\xi(\cdot)$ associated with an initial measure $\xi$. For this, given a locally bounded, Borel measurable function $g: \mathbb{R}_+ \longrightarrow \mathbb{R}_+$ and a right-continuous function $U: \mathbb{R}_+ \longrightarrow \mathbb{R}_+$ that is locally of bounded variation, let

$$(g * U)(u) = \int_{[0,u]} g(u-s)\, dU(s) \qquad \text{for all } u \geq 0.$$

Note that, by convention, the contribution to the above integral is $g(u)U(0)$ at $s = 0$ whenever $U(0) \neq 0$. Let

$$U_{\mathrm{e}}(u) = \sum_{i=0}^\infty (F_{\mathrm{e}}^{*i})(u) \qquad \text{for all } u \geq 0,$$

where $F_{\mathrm{e}}^{*0}(\cdot) \equiv 1$ and $F_{\mathrm{e}}^{*i}(\cdot) = (F_{\mathrm{e}}^{*(i-1)} * F_{\mathrm{e}})(\cdot)$ for each $i \in \{1, 2, \ldots\}$. For $\xi \in \mathcal{M}_{\mathrm{F}}^{c,p}$, define

$$(2.4) \qquad H_\xi(x) = \int_0^x \langle \mathbb{1}_{(y,\infty)}, \xi \rangle\, dy \qquad \text{for all } x \in \mathbb{R}_+.$$

It is easily verified that, for each $x \in \mathbb{R}_+$, $H_\xi(x) = \langle \chi \wedge x, \xi \rangle$, which explains why $H_\xi$ is referred to as the truncated initial workload function. Since $\xi$ has no atoms, the integrand in (2.4) is continuous. Thus, for $\xi \in \mathcal{M}_{\mathrm{F}}^{c,p}$, $H_\xi(\cdot)$ is continuously differentiable with

$$(2.5) \qquad H'_\xi(x) = \langle \mathbb{1}_{(x,\infty)}, \xi \rangle \qquad \text{for all } x \in \mathbb{R}_+.$$

In Lemma 4.4 of [10], it was shown that $\bar{S}(\cdot)$ maps $[0, \infty)$ onto $[0, \infty)$. Since $\bar{S}(\cdot)$ is also continuously differentiable and strictly increasing, it has a functional inverse, with the same properties as $\bar{S}(\cdot)$, defined on $[0, \infty)$ by

$$\bar{T}(u) = \bar{S}^{-1}(u) = \inf\{t \geq 0 : \bar{S}(t) > u\} \qquad \text{for all } u \geq 0.$$

If we let $x$ tend to $\infty$ in (2.3), execute the time change $u = \bar{S}(t)$ and then use (2.2), we obtain a convolution equation for $\bar{T}'(\cdot)$. The solution of this is

$$\bar{T}'(u) = (H'_\xi * U_\mathrm{e})(u) \qquad \text{for all } u \geq 0, \tag{2.6}$$

from which it follows that

$$\bar{T}(u) = (H_\xi * U_\mathrm{e})(u) \qquad \text{for all } u \geq 0. \tag{2.7}$$

For the full details of this derivation, see Lemma 4.4 of [10]. The convolution representation (2.7) is key to many of the developments in this paper.

In the next lemma, we use the fact that $\bar{T}(\cdot) = \bar{S}^{-1}(\cdot)$ to express $\bar{Z}(\cdot)$ as a time change of $\bar{T}'(\cdot)$ and to express (2.3) as a time change of a renewal equation.

LEMMA 2.1. *Let $\xi \in \mathcal{M}_\mathrm{F}^{c,p}$. Then, for each $t \geq 0$,*

$$\bar{Z}(t) = (H'_\xi * U_\mathrm{e})(\bar{S}(t)), \tag{2.8}$$

*and, for each $t \geq 0$ and $x \in \mathbb{R}_+$,*

$$\langle \mathbb{1}_{[0,x]}, \bar{\mu}_\xi(t) \rangle = \langle \mathbb{1}_{(\bar{S}(t), \bar{S}(t)+x]}, \xi \rangle + ((G^x * H_\xi) * U_\mathrm{e})(\bar{S}(t)). \tag{2.9}$$

PROOF. To verify (2.8), use the fact that $\bar{T}(\cdot) = \bar{S}^{-1}(\cdot)$ together with (2.2) to obtain, for each $t \geq 0$,

$$\bar{Z}(t) = \frac{1}{\bar{S}'(t)} = \bar{T}'(\bar{S}(t)).$$

This together with (2.6) implies (2.8). To verify (2.9), use the change of variables $y = \bar{S}(s)$ and the fact that $\bar{T}(\cdot) = \bar{S}^{-1}(\cdot)$ to obtain the following: for each $t \geq 0$ and $x \in \mathbb{R}_+$,

$$\int_0^t G^x(\bar{S}(t) - \bar{S}(s))\,ds = \int_0^{\bar{S}(t)} G^x(\bar{S}(t) - y)\,d\bar{T}(y) = (G^x * \bar{T})(\bar{S}(t)). \tag{2.10}$$

Substituting (2.10) into (2.3) and then using (2.7) and the associativity of the convolution operation completes the proof. □

In the next lemma, we show that, under appropriate conditions, $\bar{S}(t)$ is bounded below by a linear function for all $t$ sufficiently large.

LEMMA 2.2. *Given $\eta > 0$, there exists a finite, positive time $T^{\nu,\eta}$ depending on $\eta$ and $\nu$ such that, for all $\xi \in \mathcal{M}_\mathrm{F}^{c,p}$ with $\langle \chi, \xi \rangle < \infty$,*

$$\bar{S}(t) \geq \frac{t}{(\beta_\mathrm{e} + \eta)\langle \chi, \xi \rangle} \qquad \text{for all } t \geq \langle \chi, \xi \rangle T^{\nu,\eta}. \tag{2.11}$$

PROOF. Fix $\eta > 0$ and $\xi \in \mathcal{M}_{\mathrm{F}}^{c,p}$ such that $\langle \chi, \xi \rangle < \infty$. By the elementary renewal theorem, $U_{\mathrm{e}}(t)/t$ converges to $\beta_{\mathrm{e}}$ as $t$ tends to $\infty$ (cf. [13], Theorem 3.3.3). We note that this holds even if $\beta_{\mathrm{e}} = 0$, that is, if $\langle \chi, \nu_{\mathrm{e}} \rangle = \infty$. Thus, there exists a finite, positive time $\tilde{T}^{\nu,\eta}$ such that

$$U_{\mathrm{e}}(t) \leq (\beta_{\mathrm{e}} + \eta)t \qquad \text{for all } t \geq \tilde{T}^{\nu,\eta}.$$

Note that $\tilde{T}^{\nu,\eta}$ does not depend on $\xi$ since $U_{\mathrm{e}}(\cdot)$ does not depend on $\xi$. Moreover, since both $H_\xi$ and $U_{\mathrm{e}}$ are nondecreasing, from (2.7) and (2.4), it follows that $\bar{T}(t) \leq H_\xi(t) U_{\mathrm{e}}(t) \leq \langle \chi, \xi \rangle U_{\mathrm{e}}(t)$ for all $t \geq 0$. Thus,

(2.12) $$\bar{T}(t) \leq (\beta_{\mathrm{e}} + \eta) \langle \chi, \xi \rangle t \qquad \text{for all } t \geq \tilde{T}^{\nu,\eta}.$$

Since $\bar{S}(\cdot) = \bar{T}^{-1}(\cdot)$, it follows that

$$\bar{S}(t) \geq \frac{t}{(\beta_{\mathrm{e}} + \eta) \langle \chi, \xi \rangle} \qquad \text{for all } t \geq \bar{T}(\tilde{T}^{\nu,\eta}).$$

By (2.12), $\bar{T}(\tilde{T}^{\nu,\eta}) \leq (\beta_{\mathrm{e}} + \eta) \langle \chi, \xi \rangle \tilde{T}^{\nu,\eta}$. Setting $T^{\nu,\eta} = (\beta_{\mathrm{e}} + \eta) \tilde{T}^{\nu,\eta}$ completes the proof. $\square$

When $\beta_{\mathrm{e}} > 0$, we can set $\eta = \beta_{\mathrm{e}}$ in Lemma 2.2 to obtain the following corollary.

COROLLARY 2.3. *If $\langle \chi^2, \nu \rangle < \infty$, then there exists a positive, finite time $T^\nu$ such that, for all $\xi \in \mathcal{M}_{\mathrm{F}}^{c,p}$ with $\langle \chi, \xi \rangle < \infty$,*

(2.13) $$\bar{S}(t) \geq \frac{t}{2\beta_{\mathrm{e}} \langle \chi, \xi \rangle} \qquad \text{for all } t \geq \langle \chi, \xi \rangle T^\nu.$$

**3. The invariant manifold.** Theorem 1.1 is proved in this section. For this, we begin with the following proposition.

PROPOSITION 3.1. *For each $g \in \mathcal{C}$,*

(3.1) $$\alpha \langle g, \nu \rangle = \langle g', \nu_{\mathrm{e}} \rangle.$$

PROOF. Fix $g \in \mathcal{C}$. Note that (3.1) may be rewritten as

(3.2) $$\alpha \int_{\mathbb{R}_+} g(x)\, dF(x) = \int_{\mathbb{R}_+} g'(x) f_{\mathrm{e}}(x)\, dx.$$

Recall that for real-valued, right-continuous functions $U(\cdot)$ and $V(\cdot)$ on $\mathbb{R}_+$, which are locally of bounded variation and such that at least one of $U$ or $V$ is continuous, we have the following integration by parts formula: for all $0 \leq a < b < \infty$,

(3.3) $$\int_{(a,b]} V(x)\, dU(x) + \int_{(a,b]} U(x)\, dV(x) = V(b)U(b) - V(a)U(a)$$

(see, e.g., [8], Theorem 3.30). To prove (3.2), use $F(0) = 0$, (3.3), $g(0) = 0$ and $g$ is bounded, together with $\lim_{y \to \infty} F(y) = 1$, to obtain

$$\begin{aligned}
\int_{\mathbb{R}_+} g(x) \, dF(x) &= \int_{(0,\infty)} g(x) \, dF(x) \\
&= \lim_{y \to \infty} \int_{(0,y]} g(x) \, dF(x) \\
&= -\lim_{y \to \infty} \int_{(0,y]} g(x) \, d(1 - F(x)) \\
&= -\lim_{y \to \infty} \left[ g(y)(1 - F(y)) - \int_{(0,y]} g'(x)(1 - F(x)) \, dx \right] \\
&= \int_{\mathbb{R}_+} g'(x)(1 - F(x)) \, dx.
\end{aligned}$$

Since $f_e(x) = \alpha(1 - F(x))$ for all $x \in \mathbb{R}_+$, (3.2) holds, and hence so does (3.1). $\square$

PROOF OF THEOREM 1.1. Recall that $\bar{\mu}_{\mathbf{0}}(\cdot) \equiv \mathbf{0}$. Thus, $\xi = \mathbf{0}$ is an invariant state. Therefore, to prove Theorem 1.1, it suffices to show that $\xi \in \mathcal{M}_F^{c,p}$ is an invariant state if and only if $\xi = c\nu_e$ for some $c \in (0,\infty)$. Suppose that $\xi \in \mathcal{M}_F^{c,p}$ is an invariant state, that is, that $\bar{\mu}_\xi(t) \equiv \xi$ for all $t \geq 0$. Fix $t > 0$. Then, since $\bar{\mu}_\xi(\cdot)$ is a fluid model solution, it follows from (1.2) that, for any function $g \in \mathcal{C}$,

(3.4)  $$\int_0^t \frac{\langle g', \bar{\mu}_\xi(s) \rangle}{\langle \mathbb{1}, \bar{\mu}_\xi(s) \rangle} \, ds = \alpha t \langle g, \nu \rangle.$$

Since $\bar{\mu}_\xi(\cdot) \equiv \xi$, it follows that, for each $0 \leq s \leq t$, the numerator and the denominator of the integrand in (3.4) are given by $\langle g', \xi \rangle$ and $\langle \mathbb{1}, \xi \rangle$, respectively. Therefore, (3.4) simplifies to

$$\langle g', \xi \rangle = \langle \mathbb{1}, \xi \rangle \alpha \langle g, \nu \rangle.$$

By (3.1), the right-hand side of the above expression is given by $\langle \mathbb{1}, \xi \rangle \langle g', \nu_e \rangle$. Thus,

(3.5)  $$\langle g', \xi \rangle = c \langle g', \nu_e \rangle \quad \text{for all } g \in \mathcal{C},$$

where $c = \langle \mathbb{1}, \xi \rangle$. It turns out that $\mathcal{C}$ is a sufficiently rich class of functions in order for (3.5) to imply that $\xi = c\nu_e$, where $c = \langle \mathbb{1}, \xi \rangle$. To see this, fix $x \in (0, \infty)$. For $\varepsilon \in (0, x/2)$, let $g_\varepsilon \in \mathbf{C}_b^1(\mathbb{R}_+)$ such that $0 \leq g_\varepsilon' \leq 1$,

$$g_\varepsilon'(y) = \begin{cases} 1, & \text{if } y \in (\varepsilon, x - \varepsilon), \\ 0, & \text{if } y \in [0, \varepsilon/2] \cup [x - \varepsilon/2, \infty), \end{cases}$$

and $g_\varepsilon(y) = \int_0^y g'_\varepsilon(z)\,dz$. Then $g_\varepsilon \in \mathcal{C}$. Therefore, by (3.5), $\langle g'_\varepsilon, \xi \rangle = c\langle g'_\varepsilon, \nu_\mathrm{e} \rangle$, where $c = \langle \mathbb{1}, \xi \rangle$. By letting $\varepsilon$ tend to 0, it follows from bounded convergence that $\langle \mathbb{1}_{(0,x)}, \xi \rangle = c\langle \mathbb{1}_{(0,x)}, \nu_\mathrm{e} \rangle$ for all $x \in (0, \infty)$. Since neither $\xi$ nor $\nu_\mathrm{e}$ has an atom at the origin, $\xi = c\nu_\mathrm{e}$, where $c = \langle \mathbb{1}, \xi \rangle$. This completes the proof of the "only if" part of the theorem.

For the proof of the "if" part of the theorem, we must show that, if $\xi = c\nu_\mathrm{e}$ for some $c \in (0, \infty)$, then $\xi$ is an invariant state, that is, that $\bar\mu_\xi(\cdot) \equiv \xi$. For this, let $\bar\mu(\cdot) \equiv \xi$, where $\xi = c\nu_\mathrm{e}$ for some $c \in (0, \infty)$. It suffices to show that $\bar\mu(\cdot)$ satisfies (1.2). Obviously, for each $g \in \mathcal{C}$,

$$(3.6) \qquad \langle g, \bar\mu(t) \rangle = \langle g, \xi \rangle = \langle g, \bar\mu(0) \rangle \qquad \text{for all } t \geq 0.$$

By the definition of $\bar\mu(\cdot)$ and (3.1), for each $g \in \mathcal{C}$,

$$(3.7) \qquad \int_0^t \frac{\langle g', \bar\mu(s) \rangle}{\langle \mathbb{1}, \bar\mu(s) \rangle}\,ds = t\langle g', \nu_\mathrm{e} \rangle = \alpha t \langle g, \nu \rangle \qquad \text{for all } t \geq 0.$$

By combining (3.6) and (3.7), we see that $\bar\mu(\cdot)$ satisfies (1.2), as desired. □

**4. Weak convergence to the invariant manifold.** Theorem 1.2 is proved in this section. For this, note that, by Lemma 2.1, for $t \geq 0$ and $x \in \mathbb{R}_+$, $\bar Z(t)$ and $\langle \mathbb{1}_{(0,x]}, \bar\mu_\xi(t) \rangle$ can be expressed in terms of convolutions involving the renewal function $U_\mathrm{e}(\cdot)$. Under suitable conditions, the key renewal theorem characterizes the asymptotic behavior of such convolutions. Since $F_\mathrm{e}$ is nonarithmetic, the key renewal theorem implies that, for any Borel measurable function $g : \mathbb{R}_+ \longrightarrow \mathbb{R}_+$ that is directly Riemann integrable (see below for the definition),

$$(4.1) \qquad \lim_{z \to \infty} (g * U_\mathrm{e})(z) = \beta_\mathrm{e} \int_0^\infty g(x)\,dx,$$

(cf. [7], Chapter 11, page 363).

To apply the key renewal theorem, we will need to verify that certain functions are directly Riemann integrable. We begin by recalling the definition of the latter and some related facts. For $g : \mathbb{R}_+ \longrightarrow \mathbb{R}_+$ and $n, k \in \{1, 2, \dots\}$, let

$$m_k^n(g) = \inf\{g(z) : z \in [(k-1)/n, k/n)\}$$

and

$$M_k^n(g) = \sup\{g(z) : z \in [(k-1)/n, k/n)\},$$

and define

$$L_n(g) = \frac{1}{n} \sum_{k=1}^\infty m_k^n(g) \quad \text{and} \quad U_n(g) = \frac{1}{n} \sum_{k=1}^\infty M_k^n(g).$$

Set
$$\underline{\sigma}(g) = \limsup_{n\to\infty} L_n(g) \quad \text{and} \quad \overline{\sigma}(g) = \liminf_{n\to\infty} U_n(g).$$

The function $g$ is said to be directly Riemann integrable if $\overline{\sigma}(g) < \infty$ and $\underline{\sigma}(g) = \overline{\sigma}(g)$. Note that, for each $n$, the supremums and infimums defining $m_k^n$ and $M_k^n$ for $k \in \{1, 2, \ldots\}$ are taken over intervals of a fixed length (not varying with $k$). If $g$ is directly Riemann integrable, then $g$ is Riemann integrable and $\overline{\sigma}(g) = \int_0^\infty g(x)\,dx$. The converse is not true in general. However, if $g : \mathbb{R}_+ \longrightarrow \mathbb{R}_+$ is Riemann integrable on $[0, x]$ for all $x \in \mathbb{R}_+$ and $U_n(g) < \infty$ for some $n \in \{1, 2, \ldots\}$, then $g$ is directly Riemann integrable (cf. [7]). In particular, if $g$ is a nonincreasing, Riemann integrable function, then $g$ is directly Riemann integrable since, for all $n \in \{1, 2, \ldots\}$,
$$U_n(g) \leq \int_0^\infty g(x)\,dx + \frac{1}{n}g(0) < \infty.$$

Also, if $g_1 : \mathbb{R}_+ \longrightarrow \mathbb{R}_+$ is Riemann integrable on $[0, x]$ for all $x \in \mathbb{R}_+$ and satisfies $g_1 \leq g_2$ for some directly Riemann integrable function $g_2 : \mathbb{R}_+ \longrightarrow \mathbb{R}_+$, then $g_1$ is also directly Riemann integrable.

THEOREM 4.1. *Let $\xi \in \mathcal{M}_F^{c,p}$. If $\langle \chi, \xi \rangle < \infty$, then, for each $x \in \mathbb{R}_+$,*
$$\lim_{t\to\infty} \langle \mathbb{1}_{[0,x]}, \bar{\mu}_\xi(t) \rangle = \beta_e \langle \chi, \xi \rangle \langle \mathbb{1}_{[0,x]}, \nu_e \rangle \quad \text{and} \quad \lim_{t\to\infty} \bar{Z}(t) = \beta_e \langle \chi, \xi \rangle.$$

PROOF. Fix $\xi \in \mathcal{M}_F^{c,p}$ such that $\langle \chi, \xi \rangle < \infty$. By Lemma 2.2, $\lim_{t\to\infty} \bar{S}(t) = \infty$. Also, the continuous function $H'_\xi(\cdot)$ is directly Riemann integrable since $H'_\xi(\cdot)$ is nonincreasing and
$$\int_0^\infty H'_\xi(x)\,dx = \langle \chi, \xi \rangle < \infty.$$

These two facts together with (2.8) and (4.1) immediately imply the stated convergence result for $\bar{Z}(\cdot)$.

It remains to prove the stated convergence result for the mass on $[0, x]$ for each $x \in \mathbb{R}_+$. For this, fix $x \in \mathbb{R}_+$. Since the total mass of $\xi$ is finite and since, by Lemma 2.2, $\lim_{t\to\infty} \bar{S}(t) = \infty$, we have

(4.2) $$\limsup_{t\to\infty} \langle \mathbb{1}_{(\bar{S}(t), \bar{S}(t)+x]}, \xi \rangle \leq \lim_{t\to\infty} \langle \mathbb{1}_{(\bar{S}(t), \infty)}, \xi \rangle = 0.$$

Thus, the first term on the right-hand side of (2.9) tends to 0 as $t$ tends to $\infty$. To see that the second term on the right-hand side of (2.9) converges to the desired limit, we appeal to the key renewal theorem. For this, let $f_e^x(y) = f_e(x+y)$ for all $y \in \mathbb{R}_+$ and consider the function $(f_e^x * H_\xi)(\cdot)$. Using the fact that $H_\xi$ is continuously differentiable, it can be shown that $(f_e^x * H_\xi)(\cdot)$ is continuous, and therefore it is Riemann integrable over $[0, y]$ for each $y \in \mathbb{R}_+$.

Thus, to show that $(f_e^x * H_\xi)(\cdot)$ is directly Riemann integrable, it suffices to show that it is bounded above by a function that is directly Riemann integrable. Using the fact that both $f_e^x(\cdot)$ and $H_\xi'(\cdot)$ are nonincreasing, we obtain the following: for each $y \in \mathbb{R}_+$,

$$\begin{aligned}(f_e^x * H_\xi)(y) &= \int_0^{y/2} f_e^x(y-z) H_\xi'(z)\, dz + \int_{y/2}^{y} f_e^x(y-z) H_\xi'(z)\, dz \\ &\leq f_e^x(y/2) \int_0^{y/2} H_\xi'(z)\, dz + H_\xi'(y/2) \int_{y/2}^{y} f_e^x(y-z)\, dz \\ &\leq f_e^x(y/2) \int_0^{\infty} H_\xi'(z)\, dz + H_\xi'(y/2)(F_e(x+y/2) - F_e(x)) \\ &\leq f_e^x(y/2) \langle \chi, \xi \rangle + H_\xi'(y/2).\end{aligned}$$

As noted above, $H_\xi'(\cdot)$ is directly Riemann integrable. Since $f_e^x(\cdot)$ is nonincreasing and Riemann integrable, $f_e^x(\cdot)$ is also directly Riemann integrable. Therefore, $(f_e^x * H_\xi)(\cdot)$ is bounded above by a function that is directly Riemann integrable, and hence $(f_e^x * H_\xi)(\cdot)$ is itself directly Riemann integrable. In particular, $(f_e^x * H_\xi)(\cdot)$ is Riemann integrable over $\mathbb{R}_+$ and

$$\begin{aligned}(4.3) \quad \int_0^{\infty} (f_e^x * H_\xi)(y)\, dy &= \int_0^{\infty} \int_0^{y} f_e^x(y-z) H_\xi'(z)\, dz\, dy \\ &= \int_0^{\infty} \int_z^{\infty} f_e^x(y-z)\, dy\, H_\xi'(z)\, dz \\ &= \int_0^{\infty} \int_0^{\infty} f_e^x(y)\, dy\, H_\xi'(z)\, dz \\ &= \langle \chi, \xi \rangle \int_0^{\infty} f_e^x(y)\, dy \\ &= \langle \chi, \xi \rangle \int_x^{\infty} f_e(y)\, dy.\end{aligned}$$

Since $(G^x * H_\xi)(\cdot) = (f_e^0 * H_\xi)(\cdot) - (f_e^x * H_\xi)(\cdot)$, it immediately follows that $(G^x * H_\xi)(\cdot)$ is directly Riemann integrable. Moreover, by (4.3),

$$(4.4) \quad \int_0^{\infty} (G^x * H_\xi)(y)\, dy = \langle \chi, \xi \rangle \int_0^x f_e(y)\, dy = \langle \chi, \xi \rangle F_e(x).$$

This together with the key renewal theorem gives

$$\lim_{z \to \infty} ((G^x * H_\xi) * U_e)(z) = \beta_e \langle \chi, \xi \rangle F_e(x).$$

Since, by Lemma 2.2, $\lim_{t \to \infty} \bar{S}(t) = \infty$, it follows that

$$(4.5) \quad \lim_{t \to \infty} ((G^x * H_\xi) * U_e)(\bar{S}(t)) = \beta_e \langle \chi, \xi \rangle F_e(x).$$

Combining (2.9), (4.2) and (4.5) completes the proof. □

PROOF OF THEOREM 1.2. If $\xi = \mathbf{0}$, we have $\bar{\mu}_\xi(\cdot) \equiv \mathbf{0}$, $\beta_e \langle \chi, \xi \rangle = 0$, and the conclusion of Theorem 1.2 holds. Now take $\xi \in \mathcal{M}_F^{c,p}$ such that $\langle \chi, \xi \rangle < \infty$. First suppose that $\langle \chi^2, \nu \rangle = \infty$. Then $\beta_e = 0$, and, by Theorem 4.1, $\bar{Z}(t) \to 0$ as $t \to \infty$. Given a continuous, bounded function $g: \mathbb{R}_+ \longrightarrow \mathbb{R}$, we have, for each $t \geq 0$,

$$|\langle g, \bar{\mu}_\xi(t) \rangle| \leq \sup_{x \in \mathbb{R}_+} |g(x)| \bar{Z}(t).$$

So, it follows that $\bar{\mu}_\xi(t)$ converges weakly to the zero measure as $t \to \infty$. Next suppose that $\langle \chi^2, \nu \rangle < \infty$, that is, that $\beta_e > 0$. Recall that $\xi \neq \mathbf{0}$ implies that $\bar{Z}(t) > 0$ for all $t \geq 0$. Therefore, for each $t \geq 0$, one can divide $\bar{\mu}_\xi(t)$ by the total mass to form a probability measure. This normalization of $\bar{\mu}_\xi(t)$ for $t \geq 0$ facilitates the use of standard results on convergence in distribution. For this, for each $t \geq 0$, define the probability distribution function

$$F(t, x) = \frac{\langle \mathbb{1}_{[0,x]}, \bar{\mu}_\xi(t) \rangle}{\bar{Z}(t)} \qquad \text{for all } x \in \mathbb{R}_+.$$

Since $\beta_e \langle \chi, \xi \rangle > 0$, Theorem 4.1 implies that, for each $x \in \mathbb{R}_+$, $F(t, x) \longrightarrow F_e(x)$ as $t \to \infty$. It follows that, for any bounded, continuous function $g: \mathbb{R}_+ \longrightarrow \mathbb{R}_+$,

$$\lim_{t \to \infty} \int_{\mathbb{R}_+} g(x) \, d_x F(t, x) = \int_{\mathbb{R}_+} g(x) \, dF_e(x),$$

(cf. [5], Chapter 2, Theorem 2.2), that is, that

$$\lim_{t \to \infty} \frac{\langle g, \bar{\mu}_\xi(t) \rangle}{\bar{Z}(t)} = \langle g, \nu_e \rangle.$$

Since $\lim_{t \to \infty} \bar{Z}(t) = \beta_e \langle \chi, \xi \rangle$, it follows that, for any bounded, continuous function $g: \mathbb{R}_+ \longrightarrow \mathbb{R}_+$,

$$\lim_{t \to \infty} \langle g, \bar{\mu}_\xi(t) \rangle = \beta_e \langle \chi, \xi \rangle \langle g, \nu_e \rangle,$$

which completes the proof. □

**5. A rate of convergence in the Prohorov metric.** In this section, we prove Theorem 1.3(i), and in the next section, we prove Theorem 1.3(ii). For this, note that the conditions of part (i) [as well as those of part (ii)] imply that $\nu$ satisfies $\langle \chi^2, \nu \rangle < \infty$, and hence $\beta_e > 0$. To ease the typography, for $\beta_e > 0$ and $\xi \in \mathcal{M}_F^c$ satisfying $\langle \chi, \xi \rangle < \infty$, we use the notation

(5.1) $$\kappa = \beta_e \langle \chi, \xi \rangle.$$

We begin with Lemma 5.1, which identifies conditions that imply the result of Theorem 1.3(i). For this, recall the definitions given by (1.3) and (1.5).

LEMMA 5.1. *Let $M, \varepsilon > 0$. Suppose that $\langle \chi^2, \nu \rangle < \infty$ and that there exist a finite constant $C \geq 1$ and a finite time $T \geq 1$ such that, for all $\xi \in \mathcal{B}_\rho^{M,\varepsilon}$, $t \geq T$ and $0 < x \leq \infty$,*

$$(5.2) \qquad |\langle \mathbb{1}_{[0,x)}, \bar{\mu}_\xi(t) \rangle - \langle \mathbb{1}_{[0,x)}, \kappa \nu_{\mathrm{e}} \rangle| \leq C t^{-\varepsilon}.$$

*Then, for all $\xi \in \mathcal{B}_\rho^{M,\varepsilon}$ and $t \geq T$,*

$$\rho(\bar{\mu}_\xi(t), \kappa \nu_{\mathrm{e}}) \leq C_\rho t^{-\varepsilon/4},$$

*where $C_\rho$ is the unique positive root of the polynomial $p(y) = y^2 - (M + 4C)y - 2C$, $y \in \mathbb{R}_+$.*

PROOF. Let $M, \varepsilon > 0$ and fix $\xi \in \mathcal{B}_\rho^{M,\varepsilon}$. To prove Lemma 5.1, it suffices to show that, for each $t \geq T$ and for all nonempty closed sets $B \subset \mathbb{R}_+$,

$$(5.3) \quad \langle \mathbb{1}_B, \bar{\mu}_\xi(t) \rangle \leq \langle \mathbb{1}_{B^{\delta_t}}, \kappa \nu_{\mathrm{e}} \rangle + \delta_t \quad \text{and} \quad \langle \mathbb{1}_B, \kappa \nu_{\mathrm{e}} \rangle \leq \langle \mathbb{1}_{B^{\delta_t}}, \bar{\mu}_\xi(t) \rangle + \delta_t,$$

where $\delta_t = C_\rho t^{-\varepsilon/4}$. To verify (5.3), we begin with a simple observation. As an immediate consequence of (5.2), the fact that $\bar{\mu}_\xi(t)$ has no atoms for each $t \geq 0$ and the fact that $\nu_e$ has no atoms, it follows that, for each $0 \leq x < y \leq \infty$,

$$(5.4) \qquad |\langle \mathbb{1}_{(x,y)}, \bar{\mu}_\xi(t) \rangle - \langle \mathbb{1}_{(x,y)}, \kappa \nu_{\mathrm{e}} \rangle| \leq 2C t^{-\varepsilon} \qquad \text{for all } t \geq T.$$

We will use (5.4) in conjunction with (5.2) to verify (5.3). For this, fix a nonempty closed set $B \subset \mathbb{R}_+$ and a finite time $t \geq T$. Note that, since $B \subset B^{\delta_t}$,

$$(5.5) \qquad \langle \mathbb{1}_B, \bar{\mu}_\xi(t) \rangle \leq \langle \mathbb{1}_{B^{\delta_t}}, \bar{\mu}_\xi(t) \rangle \quad \text{and} \quad \langle \mathbb{1}_B, \kappa \nu_{\mathrm{e}} \rangle \leq \langle \mathbb{1}_{B^{\delta_t}}, \kappa \nu_{\mathrm{e}} \rangle.$$

To use (5.2) and (5.4), we will need to write $B^{\delta_t}$ as a union of intervals that are relatively open in $\mathbb{R}_+$. Since $B^{\delta_t}$ is relatively open in $\mathbb{R}_+$, it is either a finite or a countable union of relatively open, disjoint intervals. Moreover, by the definition of $B^{\delta_t}$, the length of each interval is at least $\delta_t$. Let $N$ denote the number of these intervals that have nonempty intersection with $[0, t^{\varepsilon/2})$. Then $N \leq \lceil t^{\varepsilon/2}/\delta_t \rceil$, where, for all $x \in \mathbb{R}$, $\lceil x \rceil$ denotes the smallest integer greater than or equal to $x$. Moreover, we can write

$$(5.6) \qquad B^{\delta_t} = I_1 \cup I_2 \cup \cdots \cup I_N \cup (B^{\delta_t} \cap (t^{\varepsilon/2}, \infty)),$$

where $I_i$, $i = 1, \ldots, N$, are relatively open, disjoint intervals in $\mathbb{R}_+$ such that $I_i \cap [0, t^{\varepsilon/2}) \neq \varnothing$ for $i = 1, \ldots, N$. Note that, for each $i = 1, \ldots, N$, either $I_i = [0, x)$ for some $0 < x \leq \infty$ or $I_i = (x, y)$ for some $0 \leq x < y \leq \infty$. Using (5.5) and (5.6), together with inequalities (5.2) and (5.4), we have

$$\langle \mathbb{1}_B, \bar{\mu}_\xi(t) \rangle \leq \sum_{i=1}^N \langle \mathbb{1}_{I_i}, \bar{\mu}_\xi(t) \rangle + \langle \mathbb{1}_{(t^{\varepsilon/2}, \infty)}, \bar{\mu}_\xi(t) \rangle$$

$$\leq \sum_{i=1}^{N} \langle \mathbb{1}_{I_i}, \kappa\nu_{\mathrm{e}} \rangle + \langle \mathbb{1}_{(t^{\varepsilon/2},\infty)}, \kappa\nu_{\mathrm{e}} \rangle + (N+1)2Ct^{-\varepsilon}$$

$$\leq \langle \mathbb{1}_{B^{\delta_t}}, \kappa\nu_{\mathrm{e}} \rangle + \langle \mathbb{1}_{(t^{\varepsilon/2},\infty)}, \kappa\nu_{\mathrm{e}} \rangle + (N+1)2Ct^{-\varepsilon}.$$

Similarly,

$$\langle \mathbb{1}_B, \kappa\nu_{\mathrm{e}} \rangle \leq \sum_{i=1}^{N} \langle \mathbb{1}_{I_i}, \kappa\nu_{\mathrm{e}} \rangle + \langle \mathbb{1}_{(t^{\varepsilon/2},\infty)}, \kappa\nu_{\mathrm{e}} \rangle$$

$$\leq \sum_{i=1}^{N} \langle \mathbb{1}_{I_i}, \bar{\mu}_\xi(t) \rangle + \langle \mathbb{1}_{(t^{\varepsilon/2},\infty)}, \kappa\nu_{\mathrm{e}} \rangle + N2Ct^{-\varepsilon}$$

$$\leq \langle \mathbb{1}_{B^{\delta_t}}, \bar{\mu}_\xi(t) \rangle + \langle \mathbb{1}_{(t^{\varepsilon/2},\infty)}, \kappa\nu_{\mathrm{e}} \rangle + N2Ct^{-\varepsilon}.$$

Since $\langle \chi^2, \nu \rangle < \infty$, it follows that $\langle \chi, \nu_{\mathrm{e}} \rangle < \infty$. Thus, by (5.1), $\langle \chi, \kappa\nu_{\mathrm{e}} \rangle = \langle \chi, \xi \rangle$. Therefore, $\langle \mathbb{1}_{(t^{\varepsilon/2},\infty)}, \kappa\nu_{\mathrm{e}} \rangle \leq t^{-\varepsilon/2} \langle \chi, \kappa\nu_{\mathrm{e}} \rangle = \langle \chi, \xi \rangle t^{-\varepsilon/2}$. This, together with the fact that $N \leq (t^{\varepsilon/2}/\delta_t) + 1$, gives

$$\langle \mathbb{1}_B, \bar{\mu}_\xi(t) \rangle \leq \langle \mathbb{1}_{B^{\delta_t}}, \kappa\nu_{\mathrm{e}} \rangle + \langle \chi, \xi \rangle t^{-\varepsilon/2} + \frac{2Ct^{-\varepsilon/2}}{\delta_t} + 4Ct^{-\varepsilon}$$

and

$$\langle \mathbb{1}_B, \kappa\nu_{\mathrm{e}} \rangle \leq \langle \mathbb{1}_{B^{\delta_t}}, \bar{\mu}_\xi(t) \rangle + \langle \chi, \xi \rangle t^{-\varepsilon/2} + \frac{2Ct^{-\varepsilon/2}}{\delta_t} + 2Ct^{-\varepsilon}.$$

Thus, to prove (5.3), it suffices to show that

$$\langle \chi, \xi \rangle t^{-\varepsilon/2} + 4Ct^{-\varepsilon} + \frac{2Ct^{-\varepsilon/4}}{C_\rho} \leq C_\rho t^{-\varepsilon/4}$$

or, equivalently, that

(5.7) $$(\langle \chi, \xi \rangle t^{-\varepsilon/4} + 4Ct^{-3\varepsilon/4})C_\rho + 2C \leq C_\rho^2.$$

Since $t \geq T \geq 1$ and $\xi \in \mathcal{B}_\rho^{M,\varepsilon}$,

$$\langle \chi, \xi \rangle t^{-\varepsilon/4} + 4Ct^{-3\varepsilon/4} \leq \langle \chi, \xi \rangle + 4C \leq M + 4C.$$

Moreover, since $C_\rho$ is a root of $p(\cdot)$, it follows that $(M+4C)C_\rho + 2C = C_\rho^2$. Therefore, (5.7) holds, which implies that (5.3) holds. $\square$

The next objective is to verify that the sufficient conditions in Lemma 5.1 hold under the conditions in part (i) of Theorem 1.3. It suffices to consider $\xi \neq \mathbf{0}$ since the left-hand side of (5.2) is 0 when $\xi = \mathbf{0}$. Note that, given $M, \varepsilon > 0$ and $\xi \in \mathcal{B}_\rho^{M,\varepsilon}$, by (2.9) and the fact that $\bar{\mu}_\xi(t)$ has no atoms for each $t \geq 0$, we have a useful representation for the first of the two terms that appear on the left-hand side of (5.2). Also, notice that, for each $t > 0$ and $x \in \mathbb{R}_+$,

(5.8) $\quad \langle \mathbb{1}_{(\bar{S}(t), \bar{S}(t)+x]}, \xi \rangle \leq \langle \mathbb{1}_{(\bar{S}(t), \infty)}, \xi \rangle \leq \left( \frac{1}{\bar{S}(t)} \right)^{1+\varepsilon} \langle \chi^{1+\varepsilon}, \xi \rangle.$

By Corollary 2.3, under appropriate conditions, $\bar{S}(t)$ is bounded below by a linear function of $t$ for all $t$ sufficiently large (with uniform control for $\xi \in \mathcal{B}_\rho^{M,\varepsilon}$). Thus, the most significant issue is to look at the difference between the last term in (2.9) and $\langle \mathbb{1}_{[0,x)}, \kappa \nu_e \rangle$ for each $x \in \mathbb{R}_+$. In the proof of Theorem 1.2, we used the key renewal theorem to show that, under appropriate conditions, for each $x \in \mathbb{R}_+$, the last term in (2.9) converges to $\kappa F_e(x) = \beta_e \int_0^\infty (G^x * H_\xi)(y) \, dy$ as $t$ tends to $\infty$ [cf. (5.1), (4.4) and (4.5)]. Thus, to verify that the sufficient conditions in Lemma 5.1 hold, we first identify conditions that yield a rate for this convergence that is uniform over $x \in \mathbb{R}_+$ and $\xi \in \mathcal{B}_\rho^{M,\varepsilon}$.

THEOREM 5.2. *Let $M, \varepsilon > 0$. Suppose that $\langle \chi^{2+\varepsilon}, \nu \rangle < \infty$ and that $R: [0, \infty) \longrightarrow \mathbb{R}_+$ is a nonincreasing function for which there exists a finite time $T \geq 2$ and a finite constant $C \geq 1$ such that*

(5.9) $\quad |U_e(t+s) - U_e(t) - \beta_e s| \leq CR(t) \quad \text{for all } s \in [0, 1], \ t \geq T.$

*Then there exists a finite constant $\hat{C} \geq 1$ such that, for all $\xi \in \mathcal{B}_\rho^{M,\varepsilon}$ and $x \in \mathbb{R}_+$,*

(5.10) $\quad \left| \beta_e \int_0^\infty (G^x * H_\xi)(y) \, dy - ((G^x * H_\xi) * U_e)(t) \right|$
$\quad \leq \hat{C}(t^{-\varepsilon} + R(t/2)) \quad \text{for all } t \geq T.$

Before proceeding with the proof of Theorem 5.2, we show how to use it in conjunction with Lemma 5.1 to prove part (i) of Theorem 1.3.

PROOF OF THEOREM 1.3(i). Fix $M, \varepsilon > 0$. Suppose that $\langle \chi^{2+\varepsilon}, \nu \rangle < \infty$. Let

$$R(t) = \begin{cases} 1, & 0 \leq t < 1, \\ t^{-\varepsilon}, & t \geq 1. \end{cases}$$

For each $t \geq 0$, let $D(t, s) = U_e(t+s) - U_e(t) - \beta_e s$ for all $s \geq 0$ and denote by $\mathrm{TV}_1(D(t, \cdot))$ the total variation of the function $D(t, \cdot)$ over the interval $[0, 1]$. Since, for each $t \geq 0$, $D(t, 0) = 0$, it follows that, for each $t \geq 0$,

(5.11) $\quad |U_e(t+s) - U_e(t) - \beta_e s| \leq \mathrm{TV}_1(D(t, \cdot)) \quad \text{for all } s \in [0, 1].$

Using coupling techniques, it is possible to obtain bounds on $\mathrm{TV}_1(D(t, \cdot))$ for $t$ sufficiently large (cf. [12]). For this, note that, since $\langle \chi^{2+\varepsilon}, \nu \rangle < \infty$, it follows that $\langle \chi^{1+\varepsilon}, \nu_e \rangle < \infty$. Therefore, (6.7)(ii) in III.6 of [12] with $G \equiv 1$, $U(\cdot) = U_e(\cdot)$, $B = 1$ and $\lambda = \beta_e$ implies that there exist a finite constant $C \geq 1$ and a finite time $T \geq 2$ such that

$$\mathrm{TV}_1(D(t, \cdot)) \leq Ct^{-\varepsilon} \qquad \text{for all } t \geq T.$$

This together with (5.11) implies that (5.9) holds. Thus, by Theorem 5.2, there exists a finite constant $\tilde{C} \geq 1$ [given by $(1 + 2^\varepsilon)$ times the constant $\hat{C}$ in (5.10)] such that, for all $\xi \in \mathcal{B}_\rho^{M,\varepsilon}$, $t \geq T$ and $x \in \mathbb{R}_+$,

$$(5.12) \qquad \left| \beta_e \int_0^\infty (G^x * H_\xi)(y)\, dy - ((G^x * H_\xi) * U_e)(t) \right| \leq \tilde{C} t^{-\varepsilon}.$$

Fix $\xi \in \mathcal{B}_\rho^{M,\varepsilon}$. Using (5.12), together with (2.9), the fact that $\bar{\mu}_\xi(t)$ has no atoms for each $t \geq 0$, the fact that $\nu_e$ has no atoms, (4.4), (5.8) and the fact that $T \geq 2$, we obtain, for each $x \in \mathbb{R}_+$ and $t \geq 0$ such that $\bar{S}(t) \geq T$,

$$
\begin{aligned}
&|\langle \mathbb{1}_{[0,x)}, \bar{\mu}_\xi(t) \rangle - \langle \mathbb{1}_{[0,x)}, \kappa \nu_e \rangle| \\
(5.13) \qquad &\leq \left( \frac{1}{\bar{S}(t)} \right)^{1+\varepsilon} \langle \chi^{1+\varepsilon}, \xi \rangle + \tilde{C} \left( \frac{1}{\bar{S}(t)} \right)^\varepsilon \\
&\leq (\langle \chi^{1+\varepsilon}, \xi \rangle + \tilde{C}) \left( \frac{1}{\bar{S}(t)} \right)^\varepsilon \leq (M + \tilde{C}) \left( \frac{1}{\bar{S}(t)} \right)^\varepsilon.
\end{aligned}
$$

By Corollary 2.3 and the fact that $\xi \in \mathcal{B}_\rho^{M,\varepsilon}$, there exists a finite, positive time $T^\nu$ (that does not depend on $\xi$) such that

$$(5.14) \qquad \bar{S}(t) \geq \frac{t}{2\beta_e M} \qquad \text{for all } t \geq MT^\nu.$$

Let $\tilde{T} = \max\{2\beta_e T, T^\nu\}$. Then, for all $t \geq M\tilde{T}$, (5.14) holds and $\bar{S}(t) \geq T$. This together with (5.13) gives, for all $x \in \mathbb{R}_+$ and $t \geq M\tilde{T}$,

$$(5.15) \qquad |\langle \mathbb{1}_{[0,x)}, \bar{\mu}_\xi(t) \rangle - \langle \mathbb{1}_{[0,x)}, \kappa \nu_e \rangle| \leq (M + \tilde{C})(2\beta_e M)^\varepsilon t^{-\varepsilon}.$$

By letting $x \to \infty$ in (5.15) and using the facts that, for each $t \geq 0$, $\bar{\mu}_\xi(t) \in \mathcal{M}_F$ and $\nu_e \in \mathcal{M}_F$, we see that (5.15) holds for $x = \infty$ for all $t \geq M\tilde{T}$. Therefore, (5.2) in Lemma 5.1 holds for the finite constant given by $(M + \tilde{C})(2\beta_e M)^\varepsilon$ and the finite time given by $M\tilde{T}$. So, Theorem 1.3(i) follows from Lemma 5.1. □

The final task of this section is to prove Theorem 5.2. For this, we first establish some basic properties of the functions $(G^x * H_\xi)(\cdot)$ for $x \in \mathbb{R}_+$.

PROPOSITION 5.3. *Let $\xi \in \mathcal{M}_F^{c,p}$. For each $x \in \mathbb{R}_+$, the following hold:*

(i) *For all $u \geq 0$, $0 \leq (G^x * H_\xi)(u) \leq (f_e * H_\xi)(u)$.*

(ii) *The function $(G^x * H_\xi)(\cdot)$ is absolutely continuous. In particular, for each $u \geq 0$,*

$$(5.16) \qquad (G^x * H_\xi)(u) = \int_0^u L_\xi^x(w)\,dw,$$

*where, for each $u \geq 0$,*

$$(5.17) \qquad L_\xi^x(u) = G^x(u)H'_\xi(0) - \int_0^u G^x(u-v)\xi(dv).$$

(iii) *The function $(G^x * H_\xi)(\cdot)$ is of bounded variation. In particular,*

$$\int_0^\infty |L_\xi^x(u)|\,du \leq 3\langle 1, \xi \rangle.$$

(iv) *If, for some $\varepsilon > 0$, $\langle \chi^{2+\varepsilon}, \nu \rangle < \infty$ and $\langle \chi^{1+\varepsilon}, \xi \rangle < \infty$, then, for each $x \in \mathbb{R}_+$, the following holds: for all $u > 0$,*

$$(5.18) \qquad (G^x * H_\xi)(u) \leq \int_u^\infty |L_\xi^x(v)|\,dv \leq K_\xi u^{-1-\varepsilon},$$

*where $K_\xi = (2^{1+\varepsilon} + 1)(\langle \chi^{1+\varepsilon}, \nu_e \rangle \langle 1, \xi \rangle + \langle \chi^{1+\varepsilon}, \xi \rangle)$.*

PROOF. Fix $\xi \in \mathcal{M}_F^{c,p}$ and $x \in \mathbb{R}_+$. Property (i) is immediate since $0 \leq G^x(y) \leq f_e(y)$ for all $y \geq 0$ and $H_\xi$ is nondecreasing.

To verify (5.16), note that, by Fubini's theorem, for each $u \geq 0$,

$$\int_0^u \int_0^w G^x(w-v)\xi(dv)\,dw = \int_0^u \int_v^u G^x(w-v)\,dw\,\xi(dv)$$
$$= \int_0^u \int_0^{u-v} G^x(w)\,dw\,\xi(dv).$$

Recall that, for $y \geq 0$, $H'_\xi(y) = \langle \mathbb{1}_{(y,\infty)}, \xi \rangle$. Thus, $dH'_\xi(v) = -\xi(dv)$. So we have

$$-\int_0^u \int_0^w G^x(w-v)\xi(dv)\,dw = \int_0^u \int_0^{u-v} G^x(w)\,dw\,dH'_\xi(v).$$

Thus, regarding $\int_0^{u-v} G^x(w)\,dw$ as a function of $v \in [0,u]$ and using the integration by parts formula (3.3), we obtain

$$-\int_0^u \int_0^w G^x(w-v)\xi(dv)\,dw$$
$$= -H'_\xi(0)\int_0^u G^x(w)\,dw + \int_0^u G^x(u-v)H'_\xi(v)\,dv.$$

Then (5.16) follows.

To prove (iii) and (iv), note that, by (ii) and the fact that $G^x(y) \leq f_e(y)$ for all $y \geq 0$, we have, for all $u \geq 0$,

$$\int_u^\infty |L_\xi^x(w)|\, dw \leq \int_u^\infty G^x(w) H_\xi'(0)\, dw + \int_u^\infty \int_0^w G^x(w-v)\xi(dv)\, dw$$

$$\leq \int_u^\infty f_e(w) H_\xi'(0)\, dw + \int_u^\infty \int_0^w f_e(w-v)\xi(dv)\, dw.$$

By interchanging the order of integration in the second term on the right-hand side, we obtain, for $u \geq 0$,

$$\int_u^\infty |L_\xi^x(w)|\, dw \leq H_\xi'(0)(1 - F_e(u)) + \int_0^\infty \int_{v \vee u}^\infty f_e(w-v)\, dw\, \xi(dv)$$

$$\leq H_\xi'(0)(1 - F_e(u)) + \int_0^{u/2} (1 - F_e(u-v))\xi(dv)$$

$$+ \int_{u/2}^\infty \xi(dv).$$

Now use the fact that $1 - F_e(\cdot)$ is nonincreasing to obtain, for $u \geq 0$,

$$(5.19) \quad \int_u^\infty |L_\xi^x(w)|\, dw \leq H_\xi'(0)(1 - F_e(u)) + (1 - F_e(u/2))H_\xi'(0) + H_\xi'(u/2).$$

To verify (iii), take $u = 0$ in (5.19). To prove (iv), let $\varepsilon > 0$ and assume that $\langle \chi^{2+\varepsilon}, \nu \rangle < \infty$ and $\langle \chi^{1+\varepsilon}, \xi \rangle < \infty$. Then we have $\langle \chi^{1+\varepsilon}, \nu_e \rangle < \infty$, and the second inequality in (iv) follows from (5.19) and the fact that, for all $t > 0$,

$$(5.20) \quad H_\xi'(t) \leq \langle \chi^{1+\varepsilon}, \xi \rangle t^{-1-\varepsilon} \quad \text{and} \quad 1 - F_e(t) \leq \langle \chi^{1+\varepsilon}, \nu_e \rangle t^{-1-\varepsilon}.$$

To prove the first inequality in (iv), note that, since $f_e(\cdot)$ and $H_\xi'(\cdot)$ are nonincreasing, it follows that, for $u > 0$,

$$(f_e * H_\xi)(u) \leq H_\xi'(0) \int_{u/2}^\infty f_e(y)\, dy + f_e(0) \int_{u/2}^\infty H_\xi'(y)\, dy,$$

where each term on the right-hand side of this inequality is finite for all $u \geq 0$. Therefore, by monotone convergence, each term on the right-hand side of this inequality tends to 0 as $u \to \infty$. This together with (i) implies that

$$\lim_{u \to \infty} (G^x * H_\xi)(u) = 0.$$

Consequently, by (ii), it follows that

$$(G^x * H_\xi)(u) = -\int_u^\infty L_\xi^x(v)\, dv,$$

from which the first inequality in (iv) follows. $\square$

To begin the proof of Theorem 5.2, we borrow an idea from the proof of Theorem 3.1 in [11]. The idea is to write $((G^x * H_\xi) * U_e)(t)$, for $x \in \mathbb{R}_+$ and $t \geq 0$, as a sum of integrals over intervals of length 1 [cf. (5.24)] and then to use integration by parts on each such integral [cf. (5.25)]. In this way, one obtains expressions involving the quantities that appear in (5.9). In the next proposition, a generic term in such a sum is rewritten using integration by parts. Following that, we give the proof of Theorem 5.2.

PROPOSITION 5.4. *Let $x \in \mathbb{R}_+$ and $\xi \in \mathcal{M}_F^{c,p}$. For $t \geq 1$ and $0 \leq s \leq t-1$,*

$$\int_{(s,s+1]} (G^x * H_\xi)(t-y) \, dU_e(y)$$

(5.21)
$$= (G^x * H_\xi)(t-s)(U_e(s+1) - U_e(s))$$
$$- \int_{[t-s-1,t-s)} (U_e(s+1) - U_e(t-y)) L_\xi^x(y) \, dy.$$

PROOF. Fix $x \in \mathbb{R}_+$ and $t \geq 1$ and $0 \leq s \leq t-1$. Using (5.16) and (3.3), followed by a change of variables, gives

$$\int_{(s,s+1]} (G^x * H_\xi)(t-y) \, dU_e(y)$$
$$= (G^x * H_\xi)(t-s-1) U_e(s+1) - (G^x * H_\xi)(t-s) U_e(s)$$
$$+ \int_{[t-s-1,t-s)} U_e(t-y) L_\xi^x(y) \, dy.$$

Adding and subtracting the term $(G^x * H_\xi)(t-s) U_e(s+1)$ gives

$$\int_{(s,s+1]} (G^x * H_\xi)(t-y) \, dU_e(y)$$
$$= (G^x * H_\xi)(t-s)(U_e(s+1) - U_e(s))$$
$$+ ((G^x * H_\xi)(t-s-1) - (G^x * H_\xi)(t-s)) U_e(s+1)$$
$$+ \int_{[t-s-1,t-s)} U_e(t-y) L_\xi^x(y) \, dy.$$

Using the fact that

$$(G^x * H_\xi)(t-s-1) - (G^x * H_\xi)(t-s) = -\int_{[t-s-1,t-s)} L_\xi^x(y) \, dy$$

and combining like terms gives the result. □

PROOF OF THEOREM 5.2. Fix $\xi \in \mathcal{B}_\rho^{M,\varepsilon}$ and $x \in \mathbb{R}_+$. For $t \geq T$, let $N_t = \lfloor t - T \rfloor$. For $t \geq T$, we have

(5.22) $((G^x * H_\xi) * U_e)(t) = \int_{[0,t]} (G^x * H_\xi)(t-y) \, dU_e(y) = I_1(t) + I_2(t),$

where, for $t \geq T$,
$$I_1(t) = \int_{[0,t-N_t]} (G^x * H_\xi)(t-y)\, dU_{\mathrm{e}}(y)$$

and
$$I_2(t) = \int_{(t-N_t,t]} (G^x * H_\xi)(t-y)\, dU_{\mathrm{e}}(y).$$

By parts (i) and (iv) of Proposition 5.3, we have, for $t \geq T$,
$$I_1(t) \leq K_\xi \int_{[0,t-N_t]} (t-y)^{-1-\varepsilon}\, dU_{\mathrm{e}}(y) \leq K_\xi N_t^{-1-\varepsilon} U_{\mathrm{e}}(t-N_t).$$

Note that, for $t \geq T$, we have $t - N_t \leq T + 1$. In addition, for $t \geq 2T + 2$, we have $N_t \geq t/2$. Thus, for $t \geq 2T + 2$,

(5.23) $$I_1(t) \leq 2^{1+\varepsilon} K_\xi U_{\mathrm{e}}(T+1) t^{-1-\varepsilon}.$$

For $I_2(\cdot)$, we have
$$I_2(t) = \sum_{i=1}^{N_t} \int_{(t-N_t+i-1,t-N_t+i]} (G^x * H_\xi)(t-y)\, dU_{\mathrm{e}}(y) \qquad \text{for all } t \geq T.$$
(5.24)

Thus, we can use (5.21) and then the change of variables $j = N_t - i + 1$ to obtain, for $t \geq T$,

(5.25) $$I_2(t) = \sum_{j=1}^{N_t} (G^x * H_\xi)(j)(U_{\mathrm{e}}(t+1-j) - U_{\mathrm{e}}(t-j))$$
$$- \sum_{j=1}^{N_t} \int_{[j-1,j)} (U_{\mathrm{e}}(t+1-j) - U_{\mathrm{e}}(t-y)) L_\xi^x(y)\, dy.$$

Since, for each summand, $j \leq N_t$, and since $t - N_t \geq T$, it follows that $t - j \geq T$. Thus, we can use (5.9) on each term in the first sum. Similarly, since in the integrand of each term in the second sum we have $y \leq N_t$, we can use (5.9) on each of these integrands. For $t \geq T$, this gives

(5.26) $$|I_2(t) - I_{21}(t)| \leq I_{22}(t) + I_{23}(t), \qquad t \geq T,$$

where, for $t \geq T$,
$$I_{21}(t) = \beta_{\mathrm{e}} \sum_{j=1}^{N_t} \left( (G^x * H_\xi)(j) - \int_{[j-1,j)} (y+1-j) L_\xi^x(y)\, dy \right),$$

$$I_{22}(t) = C \sum_{j=1}^{N_t} (G^x * H_\xi)(j) R(t-j),$$

$$I_{23}(t) = C \int_{[0,N_t]} R(t-y) |L_\xi^x(y)|\, dy.$$

The above representation of $I_{21}(\cdot)$ can be simplified. In fact, for $t \geq T$,

$$I_{21}(t) = \beta_{\mathrm{e}} \int_{[0,N_t)} (G^x * H_\xi)(z) \, dz. \tag{5.27}$$

To see this, write $y + 1 - j$ as $\int_{j-1}^{y} dz$ and then interchange the order of integration to obtain the following: for all $t \geq T$,

$$I_{21}(t) = \beta_{\mathrm{e}} \sum_{j=1}^{N_t} \left( (G^x * H_\xi)(j) - \int_{[j-1,j)} ((G^x * H_\xi)(j) - (G^x * H_\xi)(z)) \, dz \right)$$

$$= \beta_{\mathrm{e}} \sum_{j=1}^{N_t} \int_{[j-1,j)} (G^x * H_\xi)(z) \, dz = \beta_{\mathrm{e}} \int_{[0,N_t)} (G^x * H_\xi)(z) \, dz.$$

Let us now summarize what has been shown. By (5.22), (5.26) and (5.27), for $t \geq T$,

$$\begin{aligned}
&\left| \beta_{\mathrm{e}} \int_0^\infty (G^x * H_\xi)(z) \, dz - ((G^x * H_\xi) * U_{\mathrm{e}})(t) \right| \\
&\qquad \leq \beta_{\mathrm{e}} \int_{[N_t,\infty)} (G^x * H_\xi)(y) \, dy + I_{22}(t) + I_{23}(t) + I_1(t).
\end{aligned} \tag{5.28}$$

We have already derived an upper bound on $I_1(\cdot)$ [cf. (5.23)]. Next, we obtain estimates on the remaining terms on the right-hand side of the above inequality.

By part (iv) of Proposition 5.3 and the fact that, for $t \geq 2T + 2$, $N_t \geq t/2$, we have, for $t \geq 2T + 2$,

$$\int_{N_t}^\infty (G^x * H_\xi)(y) \, dy \leq K_\xi \int_{N_t}^\infty y^{-1-\varepsilon} \, dy = \frac{K_\xi}{\varepsilon} N_t^{-\varepsilon} \leq 2^\varepsilon \frac{K_\xi}{\varepsilon} t^{-\varepsilon}. \tag{5.29}$$

To obtain a bound on $I_{22}(\cdot)$, note that, for $t \geq T$,

$$\sum_{j=1}^{N_t} R(t-j)(G^x * H_\xi)(j)$$

$$= \sum_{j \in [1,t/2)} R(t-j)(G^x * H_\xi)(j) + \sum_{j \in [t/2,N_t]} R(t-j)(G^x * H_\xi)(j).$$

By part (iv) of Proposition 5.3, the fact that $R(\cdot)$ is nonincreasing and the fact that $t - N_t \geq T$ for $t \geq T$, we obtain, for $t \geq T$,

$$\sum_{j \in [t/2,N_t]} R(t-j)(G^x * H_\xi)(j)$$

$$\leq R(t-N_t) K_\xi \sum_{j \in [t/2,N_t]} j^{-1-\varepsilon}$$

$$\leq R(T)K_\xi \int_{t/2-1}^{\infty} y^{-1-\varepsilon}\, dy$$

$$\leq \frac{R(T)K_\xi}{\varepsilon}\left(\frac{t}{2}-1\right)^{-\varepsilon}.$$

Note that, for $t \geq 2T+2$, we have $1 - 2/t \geq T/(T+1)$, and so

$$\left(\frac{t}{2}-1\right)^{-\varepsilon} \leq \left(\frac{T+1}{T}\right)^{\varepsilon}\left(\frac{t}{2}\right)^{-\varepsilon}.$$

Thus, for $t \geq 2T+2$,

$$\sum_{j \in [t/2, N_t]} R(t-j)(G^x * H_\xi)(j) \leq \frac{R(T)K_\xi 2^\varepsilon}{\varepsilon}\left(\frac{T+1}{T}\right)^{\varepsilon} t^{-\varepsilon}.$$

Since $R(\cdot)$ is nonincreasing, it also follows from part (iv) of Proposition 5.3 that, for $t \geq 2T+2$,

$$\sum_{j \in [1, t/2)} R(t-j)(G^x * H_\xi)(j) \leq R(t/2) \sum_{j \in [1, t/2)}(G^x * H_\xi)(j)$$

$$\leq K_\xi R(t/2) \sum_{j \in [1,\infty)} j^{-1-\varepsilon}.$$

Thus, we have, for $t \geq 2T+2$,

$$(5.30) \quad I_{22}(t) \leq C\frac{R(T)K_\xi 2^\varepsilon}{\varepsilon}\left(\frac{T+1}{T}\right)^{\varepsilon} t^{-\varepsilon} + CK_\xi \left(\sum_{j=1}^{\infty} j^{-1-\varepsilon}\right) R\left(\frac{t}{2}\right).$$

We now bound $I_{23}(\cdot)$. In a similar manner to that above, for $t \geq 2T+2$,

$$I_{23}(t) = C \int_{[0, N_t)} R(t-y)|L_\xi^x(y)|\, dy$$

$$= C \int_{[0, t/2)} R(t-y)|L_\xi^x(y)|\, dy + C \int_{[t/2, N_t)} R(t-y)|L_\xi^x(y)|\, dy.$$

Then, since $R(\cdot)$ is nonincreasing, it follows from parts (iii) and (iv) of Proposition 5.3 that, for $t \geq 2T+2$,

$$I_{23}(t) \leq CR(t/2)\int_0^\infty |L_\xi^x(y)|\, dy + CR(0)\int_{[t/2,\infty)} |L_\xi^x(y)|\, dy$$

(5.31)
$$\leq 3C\langle \mathbb{1}, \xi\rangle R(t/2) + 2^{1+\varepsilon}CR(0)K_\xi t^{-1-\varepsilon}.$$

Combining (5.28)–(5.31) and (5.23) with the fact that, since $\xi \in \mathcal{B}_\rho^{M,\varepsilon}$,

$$K_\xi \leq M(2^{1+\varepsilon}+1)(\langle \chi^{1+\varepsilon}, \nu_e\rangle + 1),$$

proves the desired result. $\square$

**6. A rate of convergence in the total variation distance.** Theorem 1.3(ii) is proved in this section. For this, fix $M, \varepsilon > 0$. Throughout this section, we assume that

$$\langle \chi^{3+\varepsilon}, \nu \rangle < \infty. \tag{6.1}$$

Note that, for $\xi = \mathbf{0}$, (1.8) holds for any positive constant $C_{\text{TV}}$ and any positive time $T_{\text{TV}}$. To see this, observe that, if $\xi = \mathbf{0}$, then $\kappa = 0$ [cf. (5.1)]. Moreover, $\bar{\mu}_{\mathbf{0}}(\cdot) \equiv \mathbf{0}$, and therefore, for all $t \geq 0$, $\|\bar{\mu}_\xi(t) - \kappa \nu_{\text{e}}\|_{\text{TV}} = 0$ if $\xi = \mathbf{0}$. So it suffices to prove (1.8) for $\xi \in \mathcal{B}_{\text{TV}}^{M,\varepsilon}$ such that $\xi \neq \mathbf{0}$.

Fix $\xi \in \mathcal{B}_{\text{TV}}^{M,\varepsilon}$ and $\xi \neq \mathbf{0}$. The first order of business is to obtain an upper bound on $\|\bar{\mu}_\xi(t) - \kappa \nu_{\text{e}}\|_{\text{TV}}$, for $t \geq 0$, that is comprised of three terms [cf. (6.6)]. Then a rate of convergence to 0 as $t$ tends to $\infty$ is obtained for each of the three terms. For this, we need to introduce some notation. For a function $g: \mathbb{R}_+ \longrightarrow \mathbb{R}$ that is locally of bounded variation, let $\text{TV}_x(g)$ denote the total variation of $g$ on $[0,x]$ for each $x \in \mathbb{R}_+$. Also, denote the total variation of $g$ by $\text{TV}(g) = \lim_{x \to \infty} \text{TV}_x(g)$. Let

$$J(t,x) = \langle \mathbb{1}_{[0,x]}, \bar{\mu}_\xi(t) \rangle - \kappa \langle \mathbb{1}_{[0,x]}, \nu_{\text{e}} \rangle \qquad \text{for all } t \geq 0, x \in \mathbb{R}_+.$$

Note that, for each $t \geq 0$, neither $\bar{\mu}_\xi(t)$ nor $\nu_{\text{e}}$ charges the origin. Therefore, for all $t \geq 0$,

$$\|\bar{\mu}_\xi(t) - \kappa \nu_{\text{e}}\|_{\text{TV}} = \text{TV}(J(t,\cdot)). \tag{6.2}$$

For each $t \geq 0$, the function $J(t,\cdot)$ is readily expressed as three distinct terms. To see this, note that, by (2.3), (2.10) and the definition of $F_{\text{e}}(\cdot)$, it follows that, for all $t \geq 0$ and $x \in \mathbb{R}_+$,

$$J(t,x) = \langle \mathbb{1}_{(\bar{S}(t), \bar{S}(t)+x]}, \xi \rangle + (G^x * \bar{T})(\bar{S}(t)) - \kappa F_{\text{e}}(x). \tag{6.3}$$

Clearly, since $\bar{S}(t)$ tends to $\infty$ as $t$ tends to $\infty$, the total variation of the first term on the right-hand side of (6.3) tends to 0 as $t$ tends to $\infty$. However, individually, the total variation of the second and third terms on the right-hand side of (6.3) fails to converge to 0. Therefore, it will be necessary to take advantage of the minus sign. This can be done by expressing the third term on the right-hand side of (6.3) as a sum of two terms (cf. Lemma 6.1), the first of which combines with the second term on the right-hand side of (6.3) to form a term whose total variation tends to 0 as $t$ tends to $\infty$. For this, it will be convenient to view the convolution in the second term on the right-hand side of (6.3) as a convolution of a function with a measure. We make the following definition. Given a signed Radon measure $\zeta$ and a bounded, Borel measurable function $g: \mathbb{R}_+ \longrightarrow \mathbb{R}_+$, let

$$(g * \zeta)(x) = \int_{[0,x]} g(x-y) \zeta(dy) \qquad \text{for all } x \in \mathbb{R}_+.$$

Let $\tau$ be the Radon measure on $\mathbb{R}_+$ such that

(6.4) $$\langle \mathbb{1}_{[0,x]}, \tau \rangle = \bar{T}(x) \quad \text{for all } x \in \mathbb{R}_+.$$

Then, for all $t \geq 0$ and $x \in \mathbb{R}_+$,

(6.5) $$J(t,x) = \langle \mathbb{1}_{(\bar{S}(t),\bar{S}(t)+x]}, \xi \rangle + (G^x * \tau)(\bar{S}(t)) - \kappa F_{\mathrm{e}}(x).$$

To obtain an upper bound on $\mathrm{TV}(J(t,\cdot))$ for each $x \in \mathbb{R}_+$ we will express $F_{\mathrm{e}}(x)$ as a convolution of $G^x(\cdot)$ with Lebesgue measure on $\mathbb{R}_+$, plus a remainder term. For this, let $\ell$ denote Lebesgue measure on $\mathbb{R}_+$.

LEMMA 6.1. *For $x \in \mathbb{R}_+$,*
$$F_{\mathrm{e}}(x) = (G^x * \ell)(\bar{S}(t)) + \langle \mathbb{1}_{(\bar{S}(t),\bar{S}(t)+x]}, \nu_{\mathrm{e}} \rangle.$$

PROOF. For each $x \in \mathbb{R}_+$,
$$F_{\mathrm{e}}(x) = \int_0^x f_{\mathrm{e}}(y)\,dy = \left( \int_0^\infty f_{\mathrm{e}}(y)\,dy - \int_x^\infty f_{\mathrm{e}}(y)\,dy \right)$$
$$= \left( \int_0^\infty f_{\mathrm{e}}(y)\,dy - \int_0^\infty f_{\mathrm{e}}(x+y)\,dy \right) = \int_0^\infty G^x(y)\,dy.$$

Splitting this integral into two pieces gives, for $x \in \mathbb{R}_+$ and $t \geq 0$,
$$F_{\mathrm{e}}(x) = \int_0^{\bar{S}(t)} G^x(y)\ell(dy) + \int_{\bar{S}(t)}^\infty G^x(y)\ell(dy)$$
$$= \int_0^{\bar{S}(t)} G^x(\bar{S}(t)-y)\ell(dy) + \int_{\bar{S}(t)}^\infty (f_{\mathrm{e}}(y) - f_{\mathrm{e}}(x+y))\ell(dy)$$
$$= (G^x * \ell)(\bar{S}(t)) + \int_{\bar{S}(t)}^{x+\bar{S}(t)} f_{\mathrm{e}}(y)\ell(dy)$$
$$= (G^x * \ell)(\bar{S}(t)) + \langle \mathbb{1}_{(\bar{S}(t),\bar{S}(t)+x]}, \nu_{\mathrm{e}} \rangle. \qquad \square$$

For each $t \geq 0$ and $x \in (0, \infty)$, let
$$A(t,x) = \langle \mathbb{1}_{(\bar{S}(t),\bar{S}(t)+x]}, \xi \rangle,$$
$$B(t,x) = (G^x * \tau)(\bar{S}(t)) - \kappa(G^x * \ell)(\bar{S}(t)),$$
$$C(t,x) = \kappa \langle \mathbb{1}_{(\bar{S}(t),\bar{S}(t)+x]}, \nu_{\mathrm{e}} \rangle.$$

Also, for each $t \geq 0$, set $A(t,0) = 0$, $B(t,0) = 0$ and $C(t,0) = 0$. Note that $A(t,\cdot)$, $B(t,\cdot)$ and $C(t,\cdot)$ are right continuous. Then, by (6.5) and Lemma 6.1, for all $t \geq 0$ and $x \in \mathbb{R}_+$,
$$J(t,x) = A(t,x) + B(t,x) - C(t,x).$$

Thus, by (6.2), for all $t \geq 0$,

(6.6) $\quad \|\bar{\mu}_\xi(t) - \kappa\nu_e\|_{\mathrm{TV}} \leq \mathrm{TV}(A(t,\cdot)) + \mathrm{TV}(B(t,\cdot)) + \mathrm{TV}(C(t,\cdot)).$

To prove part (ii) of Theorem 1.3, we will bound each term on the right-hand side of (6.6) from above.

Since, for each $t \geq 0$, the functions $A(t,\cdot)$ and $C(t,\cdot)$ are nondecreasing,

(6.7) $\quad \mathrm{TV}(A(t,\cdot)) = \lim_{x \to \infty} A(t,x) - A(t,0) = \langle \mathbb{1}_{(\bar{S}(t),\infty)}, \xi \rangle,$

(6.8) $\quad \mathrm{TV}(C(t,\cdot)) = \lim_{x \to \infty} C(t,x) - C(t,0) = \kappa \langle \mathbb{1}_{(\bar{S}(t),\infty)}, \nu_e \rangle.$

By (6.7), (6.8), Chebyshev's inequality, (6.1), (2.13), (5.1) and the fact that $\xi \in \mathcal{B}_{\mathrm{TV}}^{M,\varepsilon}$, it follows that, for $t \geq MT^\nu$,

(6.9) $\quad \mathrm{TV}(A(t,\cdot)) \leq (\bar{S}(t))^{-2-\varepsilon} \langle \chi^{2+\varepsilon}, \xi \rangle \leq (2\beta_e)^{2+\varepsilon} M^{3+\varepsilon} t^{-2-\varepsilon},$

(6.10) $\quad \mathrm{TV}(C(t,\cdot)) \leq \kappa(\bar{S}(t))^{-2-\varepsilon} \langle \chi^{2+\varepsilon}, \nu_e \rangle \leq 2^{2+\varepsilon} (\beta_e M)^{3+\varepsilon} \langle \chi^{2+\varepsilon}, \nu_e \rangle t^{-2-\varepsilon}.$

Thus, for a proof of (1.8), we have obtained suitable upper bounds on the first and last terms on the right-hand side of (6.6).

The remaining task is to bound $\mathrm{TV}(B(t,\cdot))$ from above for $t$ sufficiently large. Observe that, for all $t \geq 0$ and $x \in \mathbb{R}_+$,

$$B(t,x) = (G^x * (\tau - \kappa\ell))(\bar{S}(t)) = ((f_e - f_e^x) * (\tau - \kappa\ell))(\bar{S}(t))$$
$$= \alpha((F^x - F) * (\tau - \kappa\ell))(\bar{S}(t))$$
$$= \alpha(F^x * (\tau - \kappa\ell))(\bar{S}(t)) - \alpha(F * (\tau - \kappa\ell))(\bar{S}(t)),$$

where, for each $x \in \mathbb{R}_+$, $F^x(y) = F(x+y)$ for all $y \in \mathbb{R}_+$. Note that $(F * (\tau - \kappa\ell))(\bar{S}(t))$ does not depend on $x$. Hence, it makes no contribution to the total variation of $B(t,\cdot)$. Thus, for each $t \geq 0$,

(6.11) $\quad \mathrm{TV}(B(t,\cdot)) = \alpha \mathrm{TV}(D(t,\cdot)),$

where $D(t,x) = (F^x * (\tau - \kappa\ell))(\bar{S}(t))$ for all $t \geq 0$ and $x \in \mathbb{R}_+$. To obtain a suitable upper bound on the total variation of $D(t,\cdot)$ for all $t$ sufficiently large, we introduce the following additional notation. For a signed Radon measure $\zeta$ on $\mathbb{R}_+$, let $|\zeta|$ denote the total variation measure of $\zeta$ and let $\zeta^+$ and $\zeta^-$ be nonnegative Radon measures such that $\zeta = \zeta^+ - \zeta^-$ and $|\zeta| = \zeta^+ + \zeta^-$.

LEMMA 6.2. *Let $\zeta$ be a signed Radon measure on $\mathbb{R}_+$. For fixed $r \geq 0$, define two functions $g(x) = (F^x * \zeta)(r)$ and $\hat{g}(x) = (F^x * |\zeta|)(r)$ for all $x \in \mathbb{R}_+$. Then $\mathrm{TV}(g) \leq \mathrm{TV}(\hat{g})$.*

CRITICAL PS FLUID MODEL 29PROOF. Fix $x \in \mathbb{R}_+$ and $h > 0$. We have
$$g(x+h) - g(x) = \int_0^r (F(x+h+r-y) - F(x+r-y))\zeta(dy).$$
Since $F$ is nondecreasing,
$$|g(x+h) - g(x)| \leq \int_0^r (F(x+h+r-y) - F(x+r-y))|\zeta|(dy)$$
$$= |\hat{g}(x+h) - \hat{g}(x)|.$$
The result follows from the definition of TV$(\cdot)$. $\square$

By Lemma 6.2, for each $t \geq 0$,
$$(6.12) \qquad \mathrm{TV}(D(t,\cdot)) \leq \mathrm{TV}(\hat{D}(t,\cdot)),$$
where $\hat{D}(t,x) = (F^x * |\tau - \kappa\ell|)(\bar{S}(t))$ for all $t \geq 0$ and $x \in \mathbb{R}_+$. Since $F(\cdot)$ is nondecreasing, $\hat{D}(t,\cdot)$ is also nondecreasing for each fixed $t \geq 0$. Therefore,
$$(6.13) \qquad \mathrm{TV}(\hat{D}(t,\cdot)) = \lim_{x\to\infty} \hat{D}(t,x) - \hat{D}(t,0) \qquad \text{for all } t \geq 0.$$
By monotone convergence,
$$\lim_{x\to\infty} \hat{D}(t,x) = \int_0^{\bar{S}(t)} |\tau - \kappa\ell|(dy) \qquad \text{for all } t \geq 0.$$
Therefore, by (6.13),
$$(6.14) \quad \mathrm{TV}(\hat{D}(t,\cdot)) = \int_0^{\bar{S}(t)} (1 - F(\bar{S}(t) - y))|\tau - \kappa\ell|(dy) \qquad \text{for all } t \geq 0.$$

When considering why (6.14) should be small when $t$ is large, one realizes that, for large values of the argument $y$, the measures $\tau$ and $\kappa\ell$ are close, while for small values of the argument $y$, the function $1 - F(\bar{S}(t) - y)$ is small. To take advantage of this, fix $\delta \in (0,1)$. Given $t \geq 0$, rewrite the above integral as two pieces:

$$(6.15) \qquad \int_0^{(1-\delta)\bar{S}(t)} (1 - F(\bar{S}(t) - y))|\tau - \kappa\ell|(dy),$$

$$(6.16) \qquad \int_{(1-\delta)\bar{S}(t)}^{\bar{S}(t)} (1 - F(\bar{S}(t) - y))|\tau - \kappa\ell|(dy).$$

We begin by analyzing (6.15). For each $t \geq 0$,
$$\int_0^{(1-\delta)\bar{S}(t)} (1 - F(\bar{S}(t) - y))|\tau - \kappa\ell|(dy)$$
$$(6.17) \qquad \leq (1 - F(\delta\bar{S}(t))) \int_0^{(1-\delta)\bar{S}(t)} |\tau - \kappa\ell|(dy)$$
$$\leq (1 - F(\delta\bar{S}(t)))\|\tau - \kappa\ell\|_{\mathrm{TV}}.$$

Using Chebyshev's inequality, (6.1), (2.13) and the fact that $\xi \in \mathcal{B}_{\mathrm{TV}}^{M,\varepsilon}$ gives, for each $t \geq MT^\nu$,

$$(1 - F(\delta\bar{S}(t))) \leq \langle \chi^{3+\varepsilon}, \nu \rangle (\delta\bar{S}(t))^{-3-\varepsilon} \leq \langle \chi^{3+\varepsilon}, \nu \rangle \left(\frac{2\beta_{\mathrm{e}}M}{\delta}\right)^{3+\varepsilon} t^{-3-\varepsilon}. \tag{6.18}$$

Combining (6.17) and (6.18), we have, for each $t \geq MT^\nu$,

$$\int_0^{(1-\delta)\bar{S}(t)} (1 - F(\bar{S}(t) - y))|\tau - \kappa\ell|(dy)$$
$$\leq \|\tau - \kappa\ell\|_{\mathrm{TV}} \langle \chi^{3+\varepsilon}, \nu \rangle \left(\frac{2\beta_{\mathrm{e}}M}{\delta}\right)^{3+\varepsilon} t^{-3-\varepsilon}. \tag{6.19}$$

We now analyze (6.16). For each $t \geq 0$,

$$\int_{(1-\delta)\bar{S}(t)}^{\bar{S}(t)} (1 - F(\bar{S}(t) - y))|\tau - \kappa\ell|(dy)$$
$$\leq \int_{(1-\delta)\bar{S}(t)}^{\bar{S}(t)} |\tau - \kappa\ell|(dy) \leq \langle \mathbb{1}_{[\bar{S}(t)(1-\delta),\infty)}, |\tau - \kappa\ell| \rangle. \tag{6.20}$$

Then, by combining (6.14), (6.19) and (6.20), we have, for each $t \geq MT^\nu$,

$$\mathrm{TV}(\hat{D}(t,\cdot)) \leq \|\tau - \kappa\ell\|_{\mathrm{TV}} \langle \chi^{3+\varepsilon}, \nu \rangle \left(\frac{2\beta_{\mathrm{e}}M}{\delta}\right)^{3+\varepsilon} t^{-3-\varepsilon}$$
$$+ \langle \mathbb{1}_{[\bar{S}(t)(1-\delta),\infty)}, |\tau - \kappa\ell| \rangle. \tag{6.21}$$

From (6.21), we see that what is needed are estimates on

$$\|\tau - \kappa\ell\|_{\mathrm{TV}} \quad \text{and} \quad \langle \mathbb{1}_{[r,\infty)}, |\tau - \kappa\ell| \rangle \tag{6.22}$$

for large $r$. Recall that $\kappa = \beta_{\mathrm{e}}\langle \chi, \xi \rangle$ and that $\xi \neq \mathbf{0}$. So, after factoring out $\langle \chi, \xi \rangle$ from each of the expressions in (6.22), it suffices to obtain estimates on

$$\|\tau/\langle \chi, \xi \rangle - \beta_{\mathrm{e}}\ell\|_{\mathrm{TV}} \quad \text{and} \quad \langle \mathbb{1}_{[r,\infty)}, |\tau/\langle \chi, \xi \rangle - \beta_{\mathrm{e}}\ell| \rangle,$$

for large $r$. We note that $\beta_{\mathrm{e}}\ell$ is a stationary renewal measure. To see this, consider a renewal process for which the interarrival distribution is determined by $\nu_{\mathrm{e}}$ and the initial delay distribution is determined by $(\nu_{\mathrm{e}})_{\mathrm{e}}$, where $(\nu_{\mathrm{e}})_{\mathrm{e}}$ is the excess lifetime probability measure associated with $\nu_{\mathrm{e}}$. Specifically, $(\nu_{\mathrm{e}})_{\mathrm{e}}$ is the Borel probability measure on $\mathbb{R}_+$ that is absolutely continuous with respect to Lebesgue measure on $\mathbb{R}_+$ and has density function

$$\beta_{\mathrm{e}}(1 - F_{\mathrm{e}}(x)) \qquad \text{for all } x \in \mathbb{R}_+.$$

Here note that, by (6.1), $\beta_e > 0$. This renewal process is stationary, and, for any Borel set $A \subset \mathbb{R}_+$, $\beta_e \langle \mathbb{1}_A, \ell \rangle$ is the expected number of arrivals that occur in the set $A$ (cf. [12], Chapter III.2, (2.1)). Also notice that $H_\xi(\cdot)/\langle \chi, \xi \rangle$ is a probability distribution function on $\mathbb{R}_+$. In fact, it has density function $H'_\xi(\cdot)/\langle \chi, \xi \rangle$, which makes it the excess lifetime distribution function for the Borel probability measure $\xi/\langle \mathbb{1}, \xi \rangle$ on $\mathbb{R}_+$ [cf. (2.5)]. Let $\xi_e$ denote the Borel probability measure on $\mathbb{R}_+$ associated with the distribution function $H_\xi(\cdot)/\langle \chi, \xi \rangle$. The observation that $H_\xi(\cdot)/\langle \chi, \xi \rangle$ is a probability distribution function, together with (6.4) and (2.7), implies that $\tau/\langle \chi, \xi \rangle$ is the renewal measure associated with the renewal process for which the interarrival distribution is determined by $\nu_e$ and the initial delay distribution is determined by $\xi_e$ (cf. [12], Chapter III.1, (1.4)(ii)). Therefore, what is needed are estimates on the rate at which the delayed renewal measure $\tau/\langle \chi, \xi \rangle$ converges to the stationary renewal measure $\beta_e \ell$.

One powerful tool that yields rates of convergence to stationarity for renewal measures is coupling (cf. [12]). In fact, under certain conditions, it is possible to couple two renewal processes with a common interarrival distribution so that the respective excess lifetimes agree forever after some random time $\varsigma$ called the coupling time. In our case, the common interarrival distribution is determined by $\nu_e$ and the initial delay distributions are determined by $\xi_e$ and $(\nu_e)_e$, respectively. In addition, the coupling time $\varsigma$ is finite a.s. due to (6.1) and the fact that $\xi \in \mathcal{B}_{\text{TV}}^{M,\varepsilon}$ (cf. [12], Section 5 of Chapter III). Furthermore, the results in [12] state that if the initial delay distributions and the interarrival distribution have finite $\gamma$th moments, then the coupling time $\varsigma$ has a finite $\gamma$th moment (cf. [12], Chapter III.6, (6.2)). Thus, by (6.1) and the fact that $\xi \in \mathcal{B}_{\text{TV}}^{M,\varepsilon}$, it follows that $\mathbf{E}[\varsigma^\gamma] < \infty$ for all $\gamma \in [0, 1+\varepsilon]$, where $\mathbf{E}$ denotes expected value. In fact, by carefully following the discussion on pages 83 and 84 in [12], which explains how to adapt the proof of Theorem 4.2 in Chapter II of [12] from the discrete-time setting to the continuous-time setting, and by carefully keeping track of the constants used in that argument, one can verify that, for $\gamma \in [1, 1+\varepsilon]$,

$$(6.23) \qquad \mathbf{E}[\varsigma^\gamma] \leq \frac{6^\gamma \langle \chi^{1+\gamma}, \xi \rangle}{(1+\gamma)\langle \chi, \xi \rangle} + \frac{\langle \chi^2, \xi \rangle}{2\langle \chi, \xi \rangle} C_1^\nu(\gamma) + C_2^\nu(\gamma),$$

where $C_1^\nu(\gamma)$ and $C_2^\nu(\gamma)$ are finite, positive constants that depend on $\nu$ and $\gamma$, but do not depend on $\xi$. In particular, since $\xi \in \mathcal{B}_{\text{TV}}^{M,\varepsilon}$,

$$(6.24) \qquad \langle \chi, \xi \rangle \mathbf{E}[\varsigma] \leq \left(3 + \frac{C_1^\nu(1)}{2} + C_2^\nu(1)\right) M,$$

$$(6.25) \qquad \langle \chi, \xi \rangle \mathbf{E}[\varsigma^{1+\varepsilon}] \leq \left(\frac{6^{1+\varepsilon}}{2+\varepsilon} + \frac{C_1^\nu(1+\varepsilon)}{2} + C_2^\nu(1+\varepsilon)\right) M.$$

Since it is more than a simple exercise to obtain (6.23) from the details included in [12], the verification of (6.23) is included as an Appendix here (cf. Section A.2).

Next we show how to use (6.24) and (6.25) to obtain bounds on $\mathrm{TV}(\hat{D}(t,\cdot))$ for $t \geq 0$. For this, recall that, by (6.21), it suffices to obtain bounds on the quantities that appear in (6.22). By carefully following the arguments on pages 84 and 85 of [12], it can be shown that, for $r \geq 1$,

$$(6.26) \qquad \|\tau - \kappa\ell\|_{\mathrm{TV}} \leq 2\langle \chi, \xi \rangle U_{\mathrm{e}}(1)(1 + \mathbf{E}(\varsigma)),$$

$$(6.27) \qquad \langle \mathbb{1}_{[r,\infty)}, |\tau - \kappa\ell| \rangle \leq 2\langle \chi, \xi \rangle U_{\mathrm{e}}(1)\mathbf{E}(\varsigma^{1+\varepsilon})r^{-\varepsilon}$$

(cf. Section A.3). Combining (6.26) and (6.27) with (6.21), (2.13) and the fact that $\xi \in \mathcal{B}_{\mathrm{TV}}^{M,\varepsilon}$ gives the following bound on $\mathrm{TV}(\hat{D}(\cdot,\cdot))$: for all $t \geq MT^\nu$,

$$(6.28) \quad \begin{aligned} \mathrm{TV}(\hat{D}(t,\cdot)) \\ &\leq 2\langle \chi, \xi \rangle U_{\mathrm{e}}(1)\bigg((1 + \mathbf{E}(\varsigma))\langle \chi^{3+\varepsilon}, \nu \rangle \bigg(\frac{2\beta_{\mathrm{e}}M}{\delta}\bigg)^{3+\varepsilon} t^{-3} \\ &\qquad\qquad\qquad + \mathbf{E}(\varsigma^{1+\varepsilon})\bigg(\frac{2\beta_{\mathrm{e}}M}{1-\delta}\bigg)^{\varepsilon}\bigg)t^{-\varepsilon}. \end{aligned}$$

Combining (6.28) with (6.24) and (6.25) provides a bound on $\mathrm{TV}(\hat{D}(t,\cdot))$ for $t \geq 0$ of the type that is needed to complete the proof of Theorem 1.3(ii).

PROOF OF THEOREM 1.3(ii). Fix $M, \varepsilon > 0$. If $\xi = \mathbf{0}$, it follows that $\|\bar{\mu}_\xi(t) - \kappa\nu_{\mathrm{e}}\|_{\mathrm{TV}} = 0$ for all $t \geq 0$. Therefore, it suffices to show that there exists a finite, positive constant $C_{\mathrm{TV}}$ and a finite, positive time $T_{\mathrm{TV}}$ such that (1.8) holds for all $\xi \in \mathcal{B}_{\mathrm{TV}}^{M,\varepsilon}$ such that $\xi \neq \mathbf{0}$. For this, combine (6.1), (6.6), (6.9)–(6.12) and (6.28). Then use (6.24) and (6.25). $\square$

## APPENDIX

In this appendix, we verify (6.23), (6.26) and (6.27), which were used in the proof of Theorem 1.3(ii). For this, fix $M, \varepsilon > 0$ and $\xi \in \mathcal{B}_{\mathrm{TV}}^{M,\varepsilon}$ such that $\xi \neq \mathbf{0}$. Throughout the Appendix, it is assumed that (6.1) holds.

The proofs of (6.23), (6.26) and (6.27) hinge on using the general coupling construction given in Section 5 in Chapter III of [12] to couple two renewal processes with a common interarrival distribution determined by $\nu_{\mathrm{e}}$ and initial delay distributions determined by $\xi_{\mathrm{e}}$ and $(\nu_{\mathrm{e}})_{\mathrm{e}}$, respectively. We refer to such renewal processes as $\xi_{\mathrm{e}}$-delay and stationary renewal processes, respectively. Given a $\xi_{\mathrm{e}}$-delay (resp. stationary) renewal process, let $N(\cdot)$ [resp. $N^{\mathrm{s}}(\cdot)$] denote the associated counting measure. Here the superscript s stands for stationary. Thus, for each Borel set $A \subset \mathbb{R}_+$,

$$\text{(A.1)} \qquad \mathbf{E}[N(A)] = \frac{\langle \mathbb{1}_A, \tau \rangle}{\langle \chi, \xi \rangle} \quad \text{and} \quad \mathbf{E}[N^{\text{s}}(A)] = \beta_{\text{e}} \langle \mathbb{1}_A, \ell \rangle,$$

where $\tau$ is defined by (6.4) and $\ell$ denotes Lebesgue measure. Also, for $n \in \{1, 2, \dots\}$, let $T_n$ (resp. $T_n^{\text{s}}$) denote the time of the $n$th arrival in the $\xi_{\text{e}}$-delay (resp. stationary) renewal process. By convention, set $T_0 = T_0^{\text{s}} = 0$. For $t \geq 0$, let

$$A(t) = \min\{t - T_n \geq 0 : n = 0, 1, 2, \dots\},$$
$$D(t) = \min\{T_n - t > 0 : n = 0, 1, 2, \dots\}.$$

At time $t$, $A(t)$ is the time that has elapsed since the most recent arrival in the $\xi_{\text{e}}$-delay renewal process, that is, the *age* of the most recent arrival. Similarly, $D(t)$ is the time that will elapse beginning from time $t$ until the next arrival in the $\xi_{\text{e}}$-delay renewal process, that is, the *delay* until the next arrival. For the stationary renewal process, the age process $A^{\text{s}}(\cdot)$ and the delay process $D^{\text{s}}(\cdot)$ are defined in an analogous fashion. The reason for referring to the renewal process with initial delay distribution determined by $(\nu_{\text{e}})_{\text{e}}$ as a stationary renewal process is that, for each $t \geq 0$, the distribution of $D^{\text{s}}(t)$ is equal to that of $D^{\text{s}}(0)$, which is determined by $(\nu_{\text{e}})_{\text{e}}$.

The coupling construction in [12] uses various properties of zero-delay renewal processes, which are renewal processes with initial delay distribution determined by $\delta_0$, where $\delta_0$ is the probability measure that puts one unit of mass at the origin. Given such a renewal process with interarrival distribution determined by $\nu_{\text{e}}$, the associated counting measure and other processes are defined in a manner analogous to that for the $\xi_{\text{e}}$-delay and stationary renewal processes, except that they are distinguished by the presence of a superscript z $[N^z(\cdot), T_\cdot^z, A^z(\cdot)$ and $D^z(\cdot)]$. Note that $N^z(\{0\}) = 1$, $T_1^z = 0$ almost surely and, for each $t \geq 0$,

$$\text{(A.2)} \qquad \mathbf{E}[N^z([0, t])] = U_{\text{e}}(t)$$

(cf. [12], Chapter III, (1.4)(i)). Properties of the distribution of $A^z(t)$, for $t$ sufficiently large, are used in determining the frequency of coupling attempts. Specifically, since $\nu_{\text{e}}$ has a density (which implies that it is "spread out"), by Lemma 5.1 in Chapter III of [12], there exist finite, positive constants $m$, $k$ and $T$ such that, for each $t \geq T$, the distribution of $A^z(t)$ has an absolutely continuous component for which the density is bounded below by $m$ on $[0, k]$ and $1 - F_{\text{e}}(k) > 0$. For the remainder of the Appendix, we fix such a triple, $(m, k, T)$. Note that these constants depend only on $\nu$, and not on $\xi$, since it is zero-delay renewal processes that are under consideration here.

We begin in Section A.1 by summarizing some important properties of the coupling construction given in [12]. Then, in Section A.2, we use these properties to derive a bound that is sufficient to imply (6.23) (cf. Theorem A.3). Finally (6.26) and (6.27) are verified in Section A.3.

**A.1. The coupling time.** For the case where the initial delay distributions are determined by $\xi_e$ and $(\nu_e)_e$, the interarrival distributions are determined by $\nu_e$ and the triple (associated with the interarrival distribution $\nu_e$) is given by $(m, k, T)$, the coupling construction in Section 5 of Chapter III of [12] yields a $\xi_e$-delay renewal process and a stationary renewal process, both defined on the same probability space, with certain additional properties, some of which we describe below. For this, we use the same notation for the interarrival times, age processes and delay processes associated with these two renewal processes as established at the beginning of the Appendix. In addition, we let $W_0 = 0$, $n_0 = 0$, $n_0^s = 0$ and for $i \in \{1, 2, 3, \dots\}$, we iteratively define

$$Z_{i-1} = \max\{D(W_{i-1}), D^s(W_{i-1})\}, \qquad W_i = W_{i-1} + Z_{i-1} + T,$$
$$n_i = \max\{n : T_n \leq W_i\}, \qquad n_i^s = \max\{n : T_n^s \leq W_i\}.$$

Finally, we let

(A.3) $$\mathcal{T} = \min\{i \geq 1 : D(W_i) = D^s(W_i)\}.$$

The coupling construction in [12] is such that

$$\mathbf{P}(\mathcal{T} < \infty) = 1 \quad \text{and} \quad D(t) = D^s(t) \qquad \text{for all } t \geq W_\mathcal{T}$$

(cf. [12], page 81). In fact, by (5.3) in Chapter III of [12], there exists $\delta \in (0, 1]$, which does not depend on $\xi$, such that

(A.4) $$\mathbf{P}(\mathcal{T} \geq i) \leq (1 - \delta)^{i-1} \qquad \text{for } i = 1, 2, \dots.$$

Using the fact that the interarrival distribution is determined by $\nu_e$, it is possible to show that (A.4) holds for $\delta = m^2(1 - F_e(k))k^2$. The *coupling time* $\varsigma$ is given by

(A.5) $$\varsigma = Z_0 + \sum_{i=1}^{\mathcal{T}}(T + Z_i) = Z_0 + \sum_{i=1}^{\infty} \mathbb{1}_{\{\mathcal{T} \geq i\}}(T + Z_i).$$

The times $W_i$, $i = 1, \dots, \mathcal{T}$, are the times at which coupling attempts were made. For $1 \leq i < \mathcal{T}$, each attempt was unsuccessful since $D(W_i) \neq D^s(W_i)$. However, the $\mathcal{T}$th such attempt was successful since $D(W_\mathcal{T}) = D^s(W_\mathcal{T})$. The coupling construction is such that, between successive coupling attempts, the coupled renewal processes satisfy a conditional independence property, which we now describe. For this, let, for $i \in \{0, 1, 2, \dots\}$,

$$\mathcal{F}_i = \sigma\{T_{n \wedge (n_i+1)}, T^s_{n \wedge (n_i^s+1)} : n = 0, 1, 2, \dots\}.$$

For fixed $i \in \{1,2,3,\dots\}$, conditioning on $\mathcal{F}_i$ allows the two renewal processes to be restarted at the arrival times $T_{n_i+1}$ and $T^{\mathrm{s}}_{n^{\mathrm{s}}_i+1}$, respectively. When the two renewal processes conditioned on $\mathcal{F}_i$ are restarted at their respective renewal arrival times, the coupling construction ensures that, on $\{\mathcal{T} > i\}$, they evolve as independent zero-delay renewal processes for $T + Z_i - D(W_i)$ and $T + Z_i - D^{\mathrm{s}}(W_i)$ units of time, respectively. This conditional independence property is important for the proofs given below.

**A.2. Bounds for moments of the coupling time.** In this section, we prove the following theorem, which implies (6.23).

THEOREM A.3. *Let $\gamma \in [1, 1+\varepsilon]$. Then $\varsigma$ has a finite $\gamma$th moment. Moreover,*

$$(\mathbf{E}[\varsigma^\gamma])^{1/\gamma} \leq (\mathbf{E}[Z_0^\gamma])^{1/\gamma} + (2^\gamma T^\gamma + 2^{2\gamma+1}C_1(T+1) + 2^{2\gamma+1}C_1\mathbf{E}[Z_0])^{1/\gamma}$$
$$+ \frac{((2^\gamma T^\gamma + 2^{2\gamma+1}C_1(T+1))(1-\delta) + 2^{2\gamma+2}C_1 C_2)^{1/\gamma}}{1-(1-\delta)^{1/\gamma}},$$

*where $T$ is as in Section A.1, $\delta \in (0,1]$ is as in (A.4) and $C_1$ and $C_2$ are finite, positive constants that do not depend on $\xi$ (but may depend on $\nu$ and $\gamma$).*

We begin by showing how to obtain (6.23) from Theorem A.3.

PROOF OF (6.23). Fix $\gamma \in [1, 1+\varepsilon]$. Recall that $T$ and $\delta$ do not depend on $\xi$. Let

$$K_1 = 2^\gamma T^\gamma + 2^{2\gamma+1}C_1(T+1),$$
$$K_2 = 2^{2\gamma+1}C_1,$$
$$K_3 = \frac{((2^\gamma T^\gamma + 2^{2\gamma+1}C_1(T+1))(1-\delta) + 2^{2\gamma+2}C_1 C_2)^{1/\gamma}}{1-(1-\delta)^{1/\gamma}}.$$

Thus, $K_1$, $K_2$ and $K_3$ do not depend on $\xi$. By Theorem A.3,

$$\mathbf{E}[\varsigma^\gamma] \leq ((\mathbf{E}[Z_0^\gamma])^{1/\gamma} + (K_1 + K_2\mathbf{E}[Z_0])^{1/\gamma} + K_3)^\gamma$$
(A.6) $$\leq (3\max\{(\mathbf{E}[Z_0^\gamma])^{1/\gamma}, (K_1 + K_2\mathbf{E}[Z_0])^{1/\gamma}, K_3\})^\gamma$$
$$\leq 3^\gamma \mathbf{E}[Z_0^\gamma] + 3^\gamma(K_1 + K_2\mathbf{E}[Z_0]) + 3^\gamma K_3^\gamma.$$

Since $Z_0 \leq D(0) + D^{\mathrm{s}}(0)$, it follows that $Z_0^\gamma \leq 2^\gamma(D(0))^\gamma + 2^\gamma(D^{\mathrm{s}}(0))^\gamma$. Therefore,

(A.7) $\quad \mathbf{E}[Z_0] \leq \mathbf{E}[D(0)] + \mathbf{E}[D^{\mathrm{s}}(0)] = \langle \chi, \xi_{\mathrm{e}} \rangle + \langle \chi, (\nu_{\mathrm{e}})_{\mathrm{e}} \rangle,$

(A.8) $\quad \mathbf{E}[Z_0^\gamma] \leq 2^\gamma \mathbf{E}[(D(0))^\gamma] + 2^\gamma \mathbf{E}[(D^{\mathrm{s}}(0))^\gamma] = 2^\gamma \langle \chi^\gamma, \xi_{\mathrm{e}} \rangle + 2^\gamma \langle \chi^\gamma, (\nu_{\mathrm{e}})_{\mathrm{e}} \rangle.$

It is easily verified that

(A.9) $$\langle \chi, \xi_{\mathrm{e}} \rangle = \frac{\langle \chi^2, \xi \rangle}{2\langle \chi, \xi \rangle} \quad \text{and} \quad \langle \chi^\gamma, \xi_{\mathrm{e}} \rangle = \frac{\langle \chi^{1+\gamma}, \xi \rangle}{(1+\gamma)\langle \chi, \xi \rangle}.$$

Combining (6.1) and (A.6)–(A.9) proves (6.23). □

The remaining task is to prove Theorem A.3. For this, we apply some of the general arguments given in [12] to the special case where the interarrival distribution is determined by $\nu_{\mathrm{e}}$ and the initial delays are given by $(\nu_{\mathrm{e}})_{\mathrm{e}}$ and $\xi_{\mathrm{e}}$, respectively. Since [12] does not indicate how the various constants that appear in the proofs depend on the initial delay $\xi_{\mathrm{e}}$, we provide enough details here to keep track of this dependence. For this, we follow the arguments on pages 83 and 84 in [12], filling in certain details and carefully keeping track of the constants and what they depend on. These general arguments exploit certain properties of zero-delay renewal processes. For our purposes, the statements in Lemma A.4 suffice. Note that, in Lemma A.4, it is the zero-delay renewal process that is being considered. Therefore, the constants $C_1$ and $C_2$ do not depend on $\xi$. However, they do depend on $\nu$ and the constant $C_1$ may also depend on $\gamma$.

LEMMA A.4. (i) *For each $\gamma \in [0, 2+\varepsilon]$, there exists a finite, positive constant $C_1$ such that, for all $t \geq 0$, $\mathbf{E}[(D^{\mathrm{z}}(t))^\gamma] \leq C_1(t+1)$.*

(ii) *There exists a finite, positive constant $C_2$ such that $\mathbf{E}[D^{\mathrm{z}}(t)] \leq C_2$ for all $t \geq 0$.*

PROOF. Fix $\gamma \in [0, 2+\varepsilon]$ and $t \geq 0$. Note that, by (6.1), $\langle \chi^\gamma, \nu_{\mathrm{e}} \rangle < \infty$. From page 84 of [12], it follows that

$$\mathbf{E}[(D^{\mathrm{z}}(t))^\gamma] \leq \langle \chi^\gamma, \nu_{\mathrm{e}} \rangle U_{\mathrm{e}}(t).$$

By (6.1), $\langle \chi^2, \nu_{\mathrm{e}} \rangle < \infty$. Therefore, by Lorden's inequality (cf. [12], Chapter III, (4.1)(ii)),

$$\mathbf{E}[(D^{\mathrm{z}}(t))^\gamma] \leq \frac{\langle \chi^\gamma, \nu_{\mathrm{e}} \rangle}{\langle \chi, \nu_{\mathrm{e}} \rangle} t + \frac{\langle \chi^\gamma, \nu_{\mathrm{e}} \rangle \langle \chi^2, \nu_{\mathrm{e}} \rangle}{\langle \chi, \nu_{\mathrm{e}} \rangle^2},$$

which proves (i). To prove (ii), note that

$$D^{\mathrm{z}}(t) = T^{\mathrm{z}}_{N^{\mathrm{z}}([0,t])+1} - t.$$

Therefore, by (A.2) and Wald's identity,

$$\mathbf{E}[D^{\mathrm{z}}(t)] = \langle \chi, \nu_{\mathrm{e}} \rangle (U_{\mathrm{e}}(t) + 1) - t.$$

Then, by Lorden's inequality,

$$\mathbf{E}[D^{\mathrm{z}}(t)] \leq t + \frac{\langle \chi^2, \nu_{\mathrm{e}} \rangle}{\langle \chi, \nu_{\mathrm{e}} \rangle} + \langle \chi, \nu_{\mathrm{e}} \rangle - t = \langle \chi, \nu_{\mathrm{e}} \rangle + \frac{\langle \chi^2, \nu_{\mathrm{e}} \rangle}{\langle \chi, \nu_{\mathrm{e}} \rangle},$$

which completes the proof. □

PROOF OF THEOREM A.3. Fix $\gamma \in [1, 1+\varepsilon]$. By (A.5) and Minkowski's inequality,

$$(\mathbf{E}[\varsigma^\gamma])^{1/\gamma} \leq (\mathbf{E}[Z_0^\gamma])^{1/\gamma} + \sum_{i=1}^\infty (\mathbf{E}[(T+Z_i)^\gamma \mathbb{1}_{\{\mathcal{T} \geq i\}}])^{1/\gamma}. \tag{A.10}$$

Since $\gamma \in [1, 1+\varepsilon]$, by (6.1), $\langle \chi^\gamma, (\nu_e)_e \rangle < \infty$. Also, since $\xi \in \mathcal{B}_{\mathrm{TV}}^{M,\varepsilon}$, $\langle \chi^\gamma, \xi_e \rangle < \infty$. Therefore, $\mathbf{E}[Z_0^\gamma] < \infty$ [cf. (A.8)]. The next objective is to bound $\mathbf{E}[(T+Z_i)^\gamma \mathbb{1}_{\{\mathcal{T} \geq i\}}]$ from above for each $i \in \{1, 2, 3, \dots\}$. For this, fix $i \in \{1, 2, 3, \dots\}$. Note that $\mathbb{1}_{\{\mathcal{T} \geq i\}} \in \mathcal{F}_{i-1}$. Moreover, by using the inequality $(x+y)^\gamma \leq 2^\gamma(x^\gamma + y^\gamma)$, for $x, y \in \mathbb{R}_+$, it follows that

$$\mathbf{E}[(T+Z_i)^\gamma | \mathcal{F}_{i-1}] \leq 2^\gamma T^\gamma + 2^\gamma \mathbf{E}[Z_i^\gamma | \mathcal{F}_{i-1}].$$

By definition, $Z_i = \max\{D(W_i), D^s(W_i)\}$. Therefore,

$$\mathbf{E}[Z_i^\gamma | \mathcal{F}_{i-1}] \leq 2^\gamma \mathbf{E}[(D(W_i))^\gamma | \mathcal{F}_{i-1}] + 2^\gamma \mathbf{E}[(D^s(W_i))^\gamma | \mathcal{F}_{i-1}].$$

By Lemma A.4(i) and the conditional independence property of the coupling construction, it follows that, on $\{\mathcal{T} \geq i\}$,

$$\begin{aligned}\mathbf{E}[Z_i^\gamma | \mathcal{F}_{i-1}] &\leq 2^\gamma C_1(T + Z_{i-1} - D(W_{i-1}) + 1) \\ &\quad + 2^\gamma C_1(T + Z_{i-1} - D^s(W_{i-1}) + 1) \\ &\leq 2^{\gamma+1} C_1(T + 1 + Z_{i-1}).\end{aligned}$$

Thus, on $\{\mathcal{T} \geq i\}$,

$$\mathbf{E}[(T+Z_i)^\gamma | \mathcal{F}_{i-1}] \leq 2^\gamma T^\gamma + 2^{2\gamma+1} C_1(T + 1 + Z_{i-1}).$$

This together with (A.4) implies that

$$\begin{aligned}\mathbf{E}[(T+Z_i)^\gamma \mathbb{1}_{\{\mathcal{T} \geq i\}}] &\leq (2^\gamma T^\gamma + 2^{2\gamma+1} C_1(T+1))(1-\delta)^{i-1} \\ &\quad + 2^{2\gamma+1} C_1 \mathbf{E}[Z_{i-1} \mathbb{1}_{\{\mathcal{T} \geq i\}}].\end{aligned}$$

If $i = 1$, we obtain

$$\begin{aligned}&(\mathbf{E}[(T+Z_1)^\gamma \mathbb{1}_{\{\mathcal{T} \geq 1\}}])^{1/\gamma} \\ &\qquad \leq (2^\gamma T^\gamma + 2^{2\gamma+1} C_1(T+1) + 2^{2\gamma+1} C_1 \mathbf{E}[Z_0])^{1/\gamma}.\end{aligned} \tag{A.11}$$

If $i \geq 2$, then making the observation that $\mathbb{1}_{\{\mathcal{T} \geq i\}} \leq \mathbb{1}_{\{\mathcal{T} \geq i-1\}}$ and conditioning on $\mathcal{F}_{i-2}$ gives

$$\begin{aligned}\mathbf{E}[(T+Z_i)^\gamma \mathbb{1}_{\{\mathcal{T} \geq i\}}] &\leq (2^\gamma T^\gamma + 2^{2\gamma+1} C_1(T+1))(1-\delta)^{i-1} \\ &\quad + 2^{2\gamma+1} C_1 \mathbf{E}[\mathbf{E}[Z_{i-1} | \mathcal{F}_{i-2}] \mathbb{1}_{\{\mathcal{T} \geq i-1\}}].\end{aligned}$$

Recall that $Z_{i-1} = \max\{D(W_{i-1}), D^s(W_{i-1})\} \leq D(W_{i-1}) + D^s(W_{i-1})$. Thus, if $i \geq 2$, Lemma A.4(ii) and the conditional independence property of the coupling construction imply that, on $\{\mathcal{T} \geq i-1\}$,

$$\mathbf{E}[Z_{i-1}|\mathcal{F}_{i-2}] \leq 2C_2.$$

Therefore, if $i \geq 2$,

$$
\begin{aligned}
(\mathbf{E}[(T+Z_i)^\gamma \mathbb{1}_{\{\mathcal{T} \geq i\}}])^{1/\gamma} \\
\leq ((2^\gamma T^\gamma + 2^{2\gamma+1} C_1(T+1))(1-\delta) + 2^{2\gamma+2} C_1 C_2)^{1/\gamma} \\
\times (1-\delta)^{(i-2)/\gamma}.
\end{aligned}
$$
(A.12)

Combining (A.10)–(A.12) completes the proof. $\square$

**A.3. Verification of (6.26) and (6.27).** Here we apply the general arguments given on pages 84 and 85 of [12] to the particular circumstances of interest here while keeping track of the constants to verify (6.26) and (6.27). For this, note that, for each $t \geq 0$,

$$\langle \mathbb{1}_{[t,t+1)}, |\tau - \kappa\ell| \rangle = \langle \chi, \xi \rangle \langle \mathbb{1}_{[t,t+1)}, |\tilde{\tau} - \beta_e \ell| \rangle,$$

where $\tilde{\tau} = \tau/\langle \chi, \xi \rangle$. Recall that $\tilde{\tau}$ and $\beta_e \ell$ are the renewal measures associated with the coupled $\xi_e$-delay and stationary renewal processes described in Section A.1. As noted in the introduction to the Appendix, for each $t \geq 0$, the Borel probability measure corresponding to the distribution of the delay $D^s(t)$ is given by $(\nu_e)_e$ for all $t \geq 0$. For each $t \geq 0$, let $\xi_e(t)$ denote the Borel probability measure corresponding to the distribution of the delay $D(t)$. For each $t \geq 0$, by restarting each process at time $t$, it follows that

(A.13) $$\langle \mathbb{1}_{[t,t+1)}, |\tilde{\tau} - \beta_e \ell| \rangle \leq U_e(1) \|\xi_e(t) - (\nu_e)_e\|_{\mathrm{TV}}$$

(cf. [12], Chapter III, (6.6)). By the coupling–mapping inequality (cf. [12], Chapter I, (2.12)), for each $t \geq 0$,

(A.14) $$\|\xi_e(t) - (\nu_e)_e\|_{\mathrm{TV}} \leq 2\mathbf{P}(\varsigma > t).$$

By (A.13) and (A.14), for each $r \geq 0$,

$$\langle \mathbb{1}_{[r,\infty)}, |\tilde{\tau} - \beta_e \ell| \rangle \leq \sum_{i=0}^{\infty} \langle \mathbb{1}_{[r+i,r+i+1)}, |\tilde{\tau} - \beta_e \ell| \rangle \leq U_e(1) \sum_{i=0}^{\infty} 2\mathbf{P}(\varsigma > r+i).$$

In the above inequality, for each $i \in \{0, 1, 2, \dots\}$, replace $\mathbf{P}(\varsigma > r+i)$ with

$$\sum_{j=i}^{\infty} \mathbf{P}(r+j < \varsigma \leq r+j+1),$$

interchange the order of summation and simplify, to obtain, for each $r \geq 0$,

$$(A.15) \quad \langle \mathbb{1}_{[r,\infty)}, |\tilde{\tau} - \beta_e \ell| \rangle \leq 2U_e(1) \sum_{j=0}^{\infty} (j+1)\mathbf{P}(r+j < \varsigma \leq r+j+1).$$

Letting $r = 0$ in (A.15) gives

$$\|\tilde{\tau} - \beta_e \ell\|_{\text{TV}} \leq 2U_e(1)(1 + \mathbf{E}[\varsigma]).$$

Multiplying this by $\langle \chi, \xi \rangle$ proves (6.26). To verify (6.27), fix $r \geq 1$. Then, from (A.15), using the fact that $r \geq 1$, it follows that

$$\langle \mathbb{1}_{[r,\infty)}, |\tilde{\tau} - \beta_e \ell| \rangle \leq 2U_e(1) \sum_{j=0}^{\infty} (r+j) P(r+j < \varsigma \leq r+j+1).$$

Note that, for each $j \geq 0$, $1 \leq r^{-\varepsilon}(r+j)^{\varepsilon}$. So it follows that

$$\langle \mathbb{1}_{[r,\infty)}, |\tilde{\tau} - \beta_e \ell| \rangle \leq 2U_e(1) \left( \sum_{j=0}^{\infty} (r+j)^{1+\varepsilon} P(r+j < \varsigma \leq r+j+1) \right) r^{-\varepsilon}$$

$$\leq 2U_e(1)\mathbf{E}[\varsigma^{1+\varepsilon}] r^{-\varepsilon}.$$

Multiplying the above inequality by $\langle \chi, \xi \rangle$ proves (6.27).

## REFERENCES

[1] BRAMSON, M. (1996). Convergence to equilibria for fluid models of FIFO queueing networks. *Queueing Systems Theory Appl.* **22** 5–45. MR1393404
[2] BRAMSON, M. (1997). Convergence to equilibria for fluid models of head-of-the-line proportional processor sharing queueing networks. *Queueing Systems Theory Appl.* **23** 1–26. MR1433762
[3] BRAMSON, M. (1998). State space collapse with applications to heavy traffic limits for multiclass queueing networks. *Queueing Systems Theory Appl.* **30** 89–148. MR1663763
[4] CHEN, H., KELLA, O. and WEISS, G. (1997). Fluid approximations for a processor-sharing queue. *Queueing Systems Theory Appl.* **27** 99–125. MR1490262
[5] DURRETT, R. T. (1996). *Probability*: *Theory and Examples*, 2nd ed. Duxbury, Belmont, CA. MR1609153
[6] ETHIER, S. N. and KURTZ, T. G. (1986). *Markov Processes*: *Characterization and Convergence.* Wiley, New York. MR838085
[7] FELLER, W. (1971). *An Introduction to Probability Theory and Its Applications* **2**, 2nd ed. Wiley, New York.
[8] FOLLAND, G. (1984). *Real Analysis*: *Modern Techniques and Their Applications.* Wiley, New York. MR767633
[9] GROMOLL, H. C. (2004). Diffusion approximation for a processor sharing queue in heavy traffic. *Ann. Appl. Probab.* To appear.
[10] GROMOLL, H. C., PUHA, A. L. and WILLIAMS, R. J. (2002). The fluid limit of a heavily loaded processor sharing queue. *Ann. Appl. Probab.* **12** 797–859. MR1925442

[11] HEATH, D., RESNICK, S. and SAMORODNITSKY, G. (1998). Heavy tails and long range dependence in on/off processes and associated fluid models. *Math. Oper. Res.* **23** 145–165. MR1606462
[12] LINDVALL, T. (1982). *Lectures on the Coupling Method*. Wiley, New York. MR1180522
[13] RESNICK, S. I. (1992). *Adventures in Stochastic Processes*. Birkhäuser, Boston. MR1181423
[14] WILLIAMS, R. J. (1998). Diffusion approximations for open multiclass queueing networks: Sufficient conditions involving state space collapse. *Queueing Systems Theory Appl.* **30** 27–88. MR1663759

DEPARTMENT OF MATHEMATICS
UNIVERSITY OF CALIFORNIA, SAN DIEGO
LA JOLLA, CALIFORNIA 92093-0112
USA
E-MAIL: williams@math.ucsd.edu

DEPARTMENT OF MATHEMATICS
CALIFORNIA STATE UNIVERSITY, SAN MARCOS
SAN MARCOS, CALIFORNIA 92096-0001
USA
E-MAIL: apuha@csusm.edu